\documentclass[a4paper,reqno]{amsart}

\usepackage{url} 

\usepackage[ps,dvips,xdvi,arc,all,knot,frame,curve,poly]{xy}
\usepackage[final]{graphicx}

\newcommand{\savexyscale}{\let\saved@everyxy=\everyxy}
\newcommand{\restorexyscale}{\let\everyxy=\saved@everyxy}
\newcommand{\xyscale}[1]{\everyxy={\POS#1,0,}}

\providecommand{\rgvertexscalei}{0.75}
\providecommand{\rgvertexscaleii}{1.3}
\providecommand{\rgvertexmarktwo}{\ensuremath\bullet}

\newcommand{\rgvertex}[2][\relax]{%
\ifx#1\relax%
        \POS*{\displaystyle\rgvertexmarkone}="CENTER",%
\else%
        \POS*+{\displaystyle #1}*\cir{}="CENTER",
\fi%
\rgvertexpoly{#2}%
}

\newcommand{\rgvertexpoly}[1]{\POS{%
\xypolygon#1"VERTEX"{~:{R:}~<{-}~>{}}
\xypolygon#1"AUXi"{%
~:{(\rgvertexscalei,0):}
~>{}
}
\xypolygon#1"AUXii"{%
~:{(\rgvertexscaleii,0):}
~>{}
}
}}

\newcommand{\missing}[2][\cdots]{\POS{
        0;"AUXi#2"**\dir{}?(.9)*{\displaystyle{#1}}%
}}
\newcommand{\loose}[1]{\POS{%
"VERTEX#1";"AUXii#1"**\dir{-},%
}}

\newcommand{\blackvertex}[2][\relax]{%
\ifx#1\relax%
        \POS*{\displaystyle\rgvertexmarktwo}="CENTER",%
\else%
        \POS*+{\displaystyle #1}*\cir{}="CENTER",
\fi%
\rgvertexpoly{#2}%
}

\usepackage{amsmath}
\usepackage{amssymb} 
\usepackage{amsxtra} 
\usepackage{amscd}   

\usepackage{amsfonts} 
\usepackage{mathrsfs} 

\def\bmath#1{\mbox{\boldmath$#1$}}

\usepackage{amsthm}

\theoremstyle{plain}

\newtheorem{thmsec}{Theorem}[section]
\newtheorem{lemma}{Lemma}[section]
\newtheorem{prop}{Proposition}[section]
\newtheorem*{prop*}{Proposition}

\theoremstyle{remark}

\newtheorem*{ack}{Acknowledgments}

\theoremstyle{definition}
\newtheorem{dfn}{Definition}[section]

\newcommand{\category}[1]{\mathsf{#1}} 
\newcommand{\operad}[1]{\mathcal{#1}}  

\newcommand{\setC}{\mathbb{C}}
\newcommand{\setN}{\mathbb{N}}

\newcommand{\setZ}{\mathbb{Z}}
\newcommand{\setR}{\mathbb{R}}

\newcommand{\inner}[3][]{( #2|#3 )_{#1}}   

\newcommand{\del}{\partial}
\newcommand{\ud}{\mathrm{d}}

\newcommand{\tp}[1]{\sp{\otimes #1}} 

\newcommand{\Perm}[1]{\mathfrak{S}_{#1}} 

\DeclareMathOperator{\Hom}{Hom} 
\DeclareMathOperator{\End}{End} 
\DeclareMathOperator{\Aut}{Aut} 

\DeclareMathOperator{\Id}{Id}

\newcommand{\onehalf}{\frac{1}{2}}

\newcommand{\oneof}[1]{\frac{1}{#1}}

\newcommand{\E}{\mathrm{e}} 

\newcommand{\card}[1]{\lvert#1\rvert}

\let\leq=\leqslant

\numberwithin{equation}{section} 

\newcommand{\catA}{\category{A}}        
\newcommand{\catB}{\category{B}}        
\newcommand{\catE}{\category{E}}        
\newcommand{\catF}{\category{F}}        
\newcommand{\catV}{\category{V}}        
\newcommand{\catVect}{\category{Vect}}    
\newcommand{\EndOp}{\underline{\text{\rm End}}} 
\newcommand{\cplx}[1]{{#1}_{\setC}} 

\DeclareMathOperator{\In}{In}
\DeclareMathOperator{\Out}{Out}
\DeclareMathOperator{\Leg}{Legs}

\newcommand{\gint}[2][]{%
  \left\langle#2\right\rangle_{#1}}
\newcommand{\avg}[1]{\gint{#1}}

\newcommand{\Fey}{\operad{F}}  %
\newcommand{\propF}{\mathscr{F}}  %

\DeclareMathOperator{\Src}{Src}
\DeclareMathOperator{\Tgt}{Tgt}
\newcommand{\catN}{\category{N}}

 \usepackage{prettyref}
 \newrefformat{appendix}{Appendix~\ref{#1}}
 \newrefformat{prop}{Proposition~\ref{#1}}
 \newrefformat{lemma}{Lemma~\ref{#1}}
 \newrefformat{xmp}{Example~\ref{#1}}
 \newrefformat{dfn}{Definition~\ref{#1}}
 \newrefformat{cor}{Corollary~\ref{#1}}
 \newrefformat{rem}{Remark~\ref{#1}}

\begin{document}

\title{Sums over graphs and integration over
discrete groupoids}
\date{}

\author{Domenico Fiorenza}
\address[Domenico Fiorenza]{%
    Dipartimento di Matematica ``Guido castelnuovo''\\
    Universit{\`a} di Roma ``La Sapienza'' \\
    P.le Aldo Moro, 2 \\
    00185 Roma \\
    Italy
    }
\email{fiorenza@mat.uniroma1.it}

\begin{abstract}
We show that sums over graphs such as appear in the theory of Feynman
diagrams can be seen as integrals over discrete groupoids. From this point
of view, basic combinatorial formulas of the theory of Feynman
diagrams can be interpreted as pull-back or push-forward formulas for
integrals over suitable groupoids.
\end{abstract}
\maketitle

\setcounter{tocdepth}{1} 
\tableofcontents

\section*{Introduction}
 The basic idea of the theory of Feynman diagrams is that
any Gaussian integral can be expanded into a sum
over suitable graphs:
\begin{equation*}\label{eq:prima-eq}
\int_V
\frac{P_{\Gamma}(v)}{\card{\Aut\Gamma}}
\E^{S(x_*;v)}\ud\mu(v)=\sum_{\Phi\in\Fey_{\Gamma}(0)}
\frac{Z_{x_*}^{}(\Phi)}{\card{\Aut\Phi}}\,.
\end{equation*}

The aim of this paper is to explain the above formula using the
language of integration over discrete groupoids. It should not be read as
an introduction to the theory of Feynman diagrams, nor to the theory of
integration over discrete groupoids. 
Rather, this groupoid-theoretic approach to Feynman diagrams should be
thought of as something similar to the group-theoretic
proof of the formula for the number of combinations of $k$
objects out of $n$. Indeed, that formula has a  purely
combinatorial nature and can be proved without group theory;
moreover, it would be quite strange to introduce groups and
group actions only to prove this formula. However, if one is
already familiar with basic concepts of group theory, a
group-theoretic explanation of the formula is probably
clearer than a purely combinatorial one. Similarly, the aim
of this paper is to show how combinatorial formulas which
are typical of the theory of Feynman diagrams can be easily
understood as pull-back or push-forward formulas for
integrals over suitable groupoids. 
\par
Since we need to recall basic facts from the theory of groupoids and of
Feynman diagrams, and these are apparently unrelated, we had some
difficulty in organizing the material of this paper. We chose to describe all the
general categorical features in the first part of the paper, and to show how
they appear in the context of Feynman diagrams in the second part. There are
other ways of organizing the same material: for instance one could prefer to see
each categorical construction acting on Feynman diagrams soon after its abstract
description. The diagram below shows how different sections
of the paper depend on each other, in order to allow each reader to find his
personal way through the paper.
\par
\bigskip
\noindent
\centerline{
\begin{xy}
,(0,0)*+[F]{\vbox{\hsize 2cm \noindent
Generalities\break on groupoids}};(0,-17)
*+[F]{\vbox{\hsize 1.5cm \noindent
Groupoid\break coverings}}**\dir{-}
?>*\dir{>}
,(0,0)*+[F]{\vbox{\hsize 2cm \noindent
Generalities\break on groupoids}};(35,0)
*+[F]{\vbox{\hsize 1.5cm \noindent
Feynman\break diagrams}}**\dir{-}
?>*\dir{>}
,(35,0)
*+[F]{\vbox{\hsize 1.5cm \noindent
Feynman\break diagrams}}
;(73,0)
*+[F]{\vbox{\hsize 1.5cm \noindent
Feynman\break algebras}}**\dir{-}
?>*\dir{>}
,(0,-17)
*+[F]{\vbox{\hsize 1.5cm \noindent
Groupoid\break coverings}};(35,-17)
*+[F]{\vbox{\hsize 3cm \noindent
Feynman diagrams\break with distinguished \break
sub-diagrams}}**\dir{-} ?>*\dir{>}
,(35,0)
*+[F]{\vbox{\hsize 1.5cm \noindent
Feynman\break diagrams}};(35,-17)
*+{\vbox{\hsize 3cm \noindent
\phantom{Feynman diagrams}\break\phantom{with distinguished}
\break \phantom{sub-diagrams}}}**\dir{-}
?>*\dir{>}
,(0,-17)
*+[F]{\vbox{\hsize 1.5cm \noindent
Groupoid\break coverings}};(0,-40)
*+[F]{\vbox{\hsize 2.2cm \noindent
Measures and \break integration on\break\phantom{m}discrete
\break\phantom{ii}groupoids}}**\dir{-} ?>*\dir{>}
,(0,-40)
*+[F]{\vbox{\hsize 2.2cm \noindent
Measures and \break integration on\break\phantom{m}discrete
\break\phantom{ii}groupoids}};,(35,-40)
*+[F]{\vbox{\hsize 1.7cm \noindent
Symmetric\break powers of \break groupoids}}**\dir{-} ?>*\dir{>}
,(73,0)
*+[F]{\vbox{\hsize 1.5cm \noindent
Feynman\break algebras}};(73,-40)
*+[F]{\vbox{\hsize 1.9cm \noindent
Expectation\break values}}**\dir{-} ?>*\dir{>}
,(35,0)
*+[F]{\vbox{\hsize 1.5cm \noindent
Feynman\break diagrams}};(44.5,-17)
*+{\vbox{\hsize 3cm \noindent
\phantom{Feynman diagrams}\break\phantom{with distinguished}
\break \phantom{sub-diagrams}}};(66.7,-40)
*+{\vbox{\hsize 1.9cm \noindent
\phantom{Expectation}\break \phantom{values}}}**\dir{-}
?>*\dir{>}
,(0,-40)
*+[F]{\vbox{\hsize 2.2cm \noindent
Measures and \break integration on\break\phantom{m}discrete
\break\phantom{ii}groupoids}};(35,-40)
*+[F]{\vbox{\hsize 1.7cm \noindent
Symmetric\break powers of \break groupoids}}**\dir{-}
?>*\dir{>}
,(35,-40)
*+[F]{\vbox{\hsize 1.7cm \noindent
Symmetric\break powers of \break groupoids}};(73,-40)
*+{\vbox{\hsize 1.9cm \noindent
\phantom{Expectation}\break \phantom{values}}}**\dir{-}
?>*\dir{>}
,(73,-40)
*+[F]{\vbox{\hsize 1.9cm \noindent
Expectation\break values}};(73,-63)
*+[F]{\vbox{\hsize 2.1cm \noindent
Free energy \break  and partition \break functions}}**\dir{-}
?>*\dir{>}
,(73,-63)
*+[F]{\vbox{\hsize 2.1cm \noindent
Free energy \break  and partition \break functions}}
;(73,-83)
*+[F]{\vbox{\hsize 3cm \noindent
Feynman diagrams \break\phantom{mi}  expansion of
\phantom{mi}\break Gaussian integrals}}**\dir{-}
?>*\dir{>}
,(0,-83)
*+[F]{\vbox{\hsize 1.5cm \noindent
Gaussian\break  integrals}};(73,-83)
*+[F]{\vbox{\hsize 3cm \noindent
Feynman diagrams \break\phantom{mi}  expansion of
\phantom{mi}\break Gaussian integrals}}**\dir{-}
?>*\dir{>}
\end{xy}}
\par
\bigskip
\par
While revising this paper, we learned of \cite{abdesselam;fdac}, whose
principal motivation is the Jacobian conjecture, who carries out work
similar to ours.

\begin{ack}
Special thanks go to Riccardo Longoni, Gabriele Mondello and Riccardo
Murri for many stimulating discussions on the subject and for their
comments on an early draft of this paper, and to the Referee whose
comments were essential to make the exposition clearer and more rigorous.
\end{ack}

\xyscale{/r7mm/:}

\section{Generalities on groupoids}
We begin by recalling some standard facts about groupoids (see
\cite{sga1} for details).
A \emph{groupoid} is a small category in which every morphism is
an isomorphism. If the group \(\Aut x\) is finite for any object \(x\)
of the groupoid \(\catA\), then $\catA$ is called 
a groupoid with finite homs. The smallest full subcategory of a groupoid $\catA$
containing all the objects isomorphic to a given object $x$ is called the
connected component of $x$.  
\par
Morphisms
between groupoids are functors between them; the category of
group\-oids will be denoted by the symbol {\sf Grpd}.
The isomorphism class of an object \(x\) of the groupoid \(\catA\)
will be denoted by the symbol \([x]\), and the set of isomorphism
classes of objects of \(\catA\) by the symbol \(\mathcal A\); any morphism
\(\pi\colon\catA\to\catB\) of groupoids induces a morphism of sets
\(\pi\colon{\mathcal A}\to{\mathcal B}\). More precisely, ``isomorphism 
classes''
is a functor ${\sf Grpd}\to{\sf Sets}$. By abuse of notation, we will denote
by the same symbol $[x]$ also the connected component of $A$ containing the
object $x$ or even the set of objects of this connected component; when not
specified, the exact meaning of the symbol \([x]\) in this paper will always
be clear from the context.

 Finally, we will say that a 
groupoid
\(\catA\) has a grading on the objects if it is given a groupoid
isomorphism
\(\catA\simeq\coprod_{n\in\setN}\catA_n\).
\par
A trivial example of a groupoid is  a group. Indeed,
any group \(G\) can be seen as a groupoid with only one object. Any finite
group is therefore a groupoid with finite homs. Note that, if \(G\) and
\(H\) are
two groups, groupoids morphisms between \(G\) and \(H\) are the same 
thing as group
homomorphisms between them.

More significant groupoids arise in geometry, see \cite{brown;topology}.
For instance,
for any topological space \(X\), the Poincar\'e groupoid \(\pi_1(X)\) is 
defined as the
groupoid having points of \(X\) as objects, and paths in \(X\) modulo
homotopy equivalence as morphisms. For any point
\(x\in X\) the automorphism group of \(x\) as an object of \(\pi_1(X)\) 
is the fundamental
group of \(X\) based at \(x\): \(\Aut_{\pi_1(X)}(x)=\pi_1(X,x)\). Any 
continuous map of
topological spaces \(X\to Y\) induces a groupoid morphism \(\pi_1(X)\to 
\pi_1(Y)\). Moreover, the connected components of the Poicar\'e
groupoid $\pi(X)$ correspond precisely to the (path) connected components of
the topological space $X$.
\par
We end this short section on groupoids by recalling that the 
\emph{product} of two
groupoids \(\catA\) and
\(\catB\) is the groupoid
\(\catA\times\catB\) defined by:
\begin{align*}
\text{\rm Ob}(\catA\times\catB)&=\text{\rm Ob}(\catA)\times\text{\rm 
Ob}(\catB)\\
\Hom_{\catA\times\catB}\bigl((x^{}_1,y^{}_1),(x^{}_2,y^{}_2)\bigr)&=
\Hom_{\catA}(x^{}_1,x^{}_2)\times\Hom_{\catB}(y^{}_1,y^{}_2)\,.
\end{align*}
If the objects of \(\catA\) and \(\catB\) are graded, then
the objects of
\(\catA\times\catB\) are graded by
\[
(\catA\times\catB)_n=\coprod_{m_1+m_2=n}(\catA_{m_1}\times\catB_{m_2})
\]
\medskip
\section{Groupoid coverings}\label{sec:grpd-cov}

In this section we recall a few facts from the theory of groupoid 
coverings (see
\cite{sga1,brown;topology}). A morphism \(\pi\colon\catA\to\catB\) of
groupoids is
called a
\emph{fibration} if, for any morphism \(y^{}_1\to y^{}_2\) in \(\catB\)  and 
any \(x^{}_1\)
with \(\pi(x^{}_1)=y{}_1\) there exist a morphism
\(x^{}_1\to x^{}_2\) in \(\catA\) such that \(\pi(x^{}_1\to
x{}_2)=\{y^{}_1\to  y^{}_2\}\). Such a morphism is
called a \emph{lifting} of \(y^{}_1\to y^{}_2\). If the lifting is unique for
any $x_1$ lying over $y_1$, then the  fibration
\(\pi\) is called a \emph{groupoid covering}. As a consequence of the
uniqueness of the lifting, if
\(\pi\colon\catA\to\catB\) is a groupoid covering then the
induced morphisms \(\pi\colon\Aut x\to\Aut \pi(x)\) are injective, for any
\(x\in \text{\rm Ob}(\catA)\). This terminology has a clear origin in
topology. In fact, if
 \(X\) is a topological space and \(E\to X\) a Serre fibration, then
\(\pi_1(E)\to \pi_1(X)\) is a groupoid fibration, and if \(\tilde X\to X\) 
is a covering, then
\(\pi_1(\tilde X)\to \pi_1(X)\) is a groupoid covering.
In the case of groups, a groupoid fibration \(\pi:G\to H\) is just a
surjective group homomorphism, and a covering \(\pi:G\to H\) is a group
 isomorphism. \par
%

 If \(y\) is an object of \(\catB\), we
denote by the symbol \(\category{y}\) the full subcategory of \(\catB\)
having  only \(y\) as an
object; the symbol \(\pi^{-1}(\category{y})\) denotes the full subcategory of 
\(\catA\) whose
objects are the objects \(x\) of \(\catA\) such that \(\pi(x)=y\). Note 
that, if \(\pi\colon\catA\to\catB\) is a covering, then also
\begin{equation*}
\pi\vert_{\pi^{-1}(\category{y})}\colon\pi^{-1}(\category{y})\to \category{y}
\end{equation*}
is a covering, for any object \(y\) of \(\catB\). If \(y\) and \(y'\) are 
two objects
in the same isomorphism class in
\(\catB\), then the subcategories \(\pi^{-1}(\category{y})\) and
\(\pi^{-1}(\category{y}')\)  of \(\catA\)
are isomorphic. Moreover, there is a canonical injection
\(\Hom_{\catB}(y,y')\hookrightarrow\text{\rm
Iso}\bigl(\pi^{-1}(\category{y}),\pi^{-1}(\category{y}')\bigr)\).
Indeed, let \(\phi\colon y\to y'\) be
an isomorphism. Since \(\pi\) is a covering, for any
\(x\in\text{\rm Ob}\left(\pi^{-1}(\category{y})\right)\) there exist a unique
lifting
\(\tilde\phi\colon x\to x'\) of \(\phi\). The map
\(\tilde\phi\colon\pi^{-1}(\category{y})\to\pi^{-1}(\category{y}')\) is an
isomorphism, and the map \(\phi\mapsto \tilde\phi\) is the claimed injection.
\par
When $y=y'$, we obtain an action of \(\Aut y\) on \(\pi^{-1}(\category{y})\).
It is  immediate
to see that this action  restricts to a
transitive action on the set
\[
[x]_y:=\text{\rm Ob}\left(\pi^{-1}(\category{y})\cap[x]\right)
\]
for any \(x\in \text{\rm Ob}\left(\pi^{-1}(\category{y})\right)\).
On the other hand, since
\(\pi(x)=y\), any
\(\phi\in\Aut y\) lifts to some \(\tilde\phi\colon
x\to x'\). Therefore,
\(\phi\) stabilizes \(x\) if and only if it lifts to some
\(\psi\colon x\to x\), i.e., \(\phi=\pi(\psi)\) for some \(\psi\in \Aut
x\). This
means that \(\text{\rm Stab} x=\pi(\Aut x)\), and we have a canonical 
isomorphism
of \(\Aut y\)-sets
\begin{equation}\label{eq:action-on-fiber}
\Aut y/\pi(\Aut x) \simeq [x]_y
\end{equation}
\par
If \(p\colon \tilde X\to X\) is a covering of (path connected)
topological spaces, the isomorphism \prettyref{eq:action-on-fiber} above is just
the well known isomorphism
\begin{equation*}
p^{-1}(x_0)\simeq \pi_1( X, x_0)/p_*\pi_1(\tilde
X,\tilde x_0)\,,
\end{equation*}
where \(\tilde x_0\) is any point in the fibre \(p^{-1}(x_0)\).
\medskip
\par

\par

A covering $\pi\colon\catA\to\catB$ is called \emph{finite} if, for any
\(y\in\text{\rm Ob}(\catB)\), the pre-image
\(\pi^{-1}(\category{y})\) consists of finitely many objects.  Since the map 
\(\pi\colon\Aut
x\to \Aut y\) is injective by definition of covering, if
\(\Aut y\) is finite the isomorphism 
\prettyref{eq:action-on-fiber} implies
\begin{equation}\label{eq:cardinality}
\frac{\card{\Aut y}}{\card{\Aut x}}=\card{[x]_y}
\end{equation}
\par
 The \emph{degree} of a finite covering $\pi\colon\catA\to\catB$ is the
function
\begin{align*}
\deg\pi\colon \text{\rm Ob}(\catB)&\to \setN\\
y&\mapsto \card{\text{\rm Ob}\left(\pi^{-1}(\category{y})\right)}
\end{align*}
It is immediate to compute
\begin{equation}\label{eq:degree-cov}
\deg\pi(y)=\sum_{[x]\in\pi^{-1}([y])}\card{\text{\rm
Ob}\left(\pi^{-1}(\category{y})\cap [x]\right)}
=\sum_{[x]\in\pi^{-1}([y])}\card{[x]_y}
\end{equation}
\par
If $y\simeq y'$, then $\pi^{-1}(\category{y})$  and $\pi^{-1}(\category{y}')$
are  isomorphic, so
the degree of a covering is constant on isomorphism classes, and defines a
map
\begin{equation*}
\deg\pi:\mathcal B\to \setN\,.
\end{equation*}
A finite covering of groupoids
\(\pi\colon \catA\to \catB\) is called \emph{homogeneous}  if its degree 
is constant
on \(\mathcal B\).
\par

\medskip
\section{Measures and integration on discrete groupoids}\label{sec:int-over-grpd}

 Groupoids can be topologized; in this paper we will be interested
only in group\-oids endowed with the
discrete topology, which will be called discrete groupoids. Clearly, any
groupoid can be topologized as a discrete groupoid.
The \emph{coarse space} of a discrete groupoid \(\catA\) is the set
\(\mathcal A\) of isomorphism classes of objects of \(\catA\),
endowed with the discrete topology.
Any morphism \(\pi\colon\catA\to\catB\) induces a map
between the coarse spaces: \(\pi\colon
\mathcal A\to\mathcal B\). The
\emph{discrete measure} \(\mu_{\mathcal A}^{}\) on the coarse space of
a discrete groupoid \(\catA\) with finite homs is defined as
\begin{equation}\label{eq:measure}
\mu_{\mathcal A}^{}([x]):=\frac{1}{\card{\Aut x}}, \qquad \forall [x]\in 
\mathcal A,
\end{equation}
where \(x\) is any object in the class \([x]\).
It is immediate from the definition above that the coarse space of
\(\catA\times\catB\) is
\(\mathcal A\times \mathcal B\) and that
\begin{equation*}
\mu_{\mathcal A\times\mathcal B}^{}=\mu_{\mathcal 
A}^{}\otimes\mu_{\mathcal B}^{}
\end{equation*}

\medskip

Let now \(V\) be a \(\setC\)-vector space or a module over some commutative
unitary
\(\setC\)-algebra \(R\). We can look at
\(V\) as a trivial groupoid, by setting
\begin{equation*}
\Hom_{V}(v_1,v_2)=\begin{cases}
\{\Id_{v_1}\}&\text{ if }v_1= v_2 \\
\emptyset&\text{ if }v_1\neq v_2
\end{cases}
\end{equation*}
Clearly, a morphism \(\catA\to V\) is just a map \(\text{\rm Ob}(\catA)\to
V\) which is constant on isomorphism classes, i.e., a map \(\mathcal A\to
V\).

We are now interested in defining integrable functions on the coarse 
space \(\mathcal A\) of
a groupoid \({\sf A}\). Since the measure we are considering is 
discrete, an
integrable function will have to vanish outside a countable subset of
\(\mathcal
A\); to simplify our treatment, we assume the coarse space \(\mathcal A\)
we are going to integrate on is countable. Moreover, since
integrals will be defined by a
limiting procedure, we will have to work with topological algebras
and 
modules (see
\cite{topological-rings} for details). The base field
\(\setC\) will always be given the Euclidean topology.

\begin{dfn}
Let \({\sf A}\simeq\coprod_n{\sf A}_n\) be a discrete
groupoid with finite homs and graded objects, such that
the coarse spaces
\(\{{\mathcal
A}_n\}\) are finite sets and let \(V\) be a topological module
over a topological commutative
\(\setC\)-algebra \(R\).  We say that a morphism \(\varphi\colon\catA\to
V\) is
\emph{integrable} if the series
\begin{equation}\label{eq:partial-sum}
\sum_{n=0}^\infty\left(\sum_{[x]\in{\mathcal
A}_n}\frac{\varphi(x)}{\card{\Aut x}}\right)
\end{equation}
is convergent. If so, we write
\begin{equation*}
\int_{\mathcal A}\varphi\,\ud \mu^{}_{\mathcal A}=\sum_{[x]\in\mathcal
A}\frac{\varphi(x)}{\card{\Aut x}}
\end{equation*}
\end{dfn}
For instance, if we consider the trivial groupoid \(\sf N\)
 on the set of natural numbers and the gradation \(\{{\sf
N}_n\}\) given by \({\sf N}_n=\{n\}\), then a morphism
\begin{align*}
a\colon {\sf N}&\to \setC\\
n&\mapsto a^{}_n
\end{align*}
is integrable if and only if the series
\(
\sum_{n=0}^\infty a^{}_n
\)
is convergent.
We now prove a criterion for integrability of a morphism
\(\varphi\colon\catA\to V\).
\begin{lemma}\label{lemma:filtrated}
Let \(V\) be a graded complete \(R\)-module, and let
\(\varphi\colon {\sf A}\to V\) be a grading-preserving morphism. Then
\(\varphi\) is integrable.
\end{lemma}
\begin{proof}
Since \(V\) is a graded complete module, \(V\) is the completion of 
\(\oplus_{n=0}^\infty
V_n\) in the product topology, where \(V_n\) is the submodule of
\(V\) consisting of elements of
degree \(n\). Therefore, a sequence converges in \(V\) if and only if 
its \(n\)-degree
components converge in \(V_n\), for every \(n\). Since the morphism
\(\varphi\) is grading-preserving
and \(\mathcal A\) contains only finitely many  degree \(n\) elements 
for any
fixed \(n\),
the
\(n\)-degree component of the sequence of the partial sums
\prettyref{eq:partial-sum} stabilizes for any fixed \(n\), hence the 
statement.
\end{proof}

If $\catA$ is a graded discrete groupoid, we denote by \(\mathscr{A}_{n}\)
the free $\setC$-vector space generated by the elements of \(\mathcal 
A\) of degree \(n\),
endowed with the Euclidean topology, and by
\(\mathscr{A}\) the direct sum of these spaces:
\begin{equation*}
\mathscr{A}:=\bigoplus_{n\in\setN}\mathscr{
A}_{n}\,.
\end{equation*}
Finally, let
\(\mathscr{A}\hskip -1 em\underline{\hskip 1 em}\hskip .1 em\) be the
completion of
\(\mathscr{A}\) with respect to the product topology.

\prettyref{lemma:filtrated} immediately implies
that the natural embedding \(j_\catA\colon\mathcal A\to
\mathscr{A}\hskip -1 em\underline{\hskip 1 em}\hskip .1 em\) is
integrable, so that we can write
\begin{equation*}
\int_{\mathcal A}j_\catA\ud \mu^{}_{\mathcal A}=\sum_{[x]\in\mathcal
A}\frac{[x]}{\card{\Aut x}}\,;
\end{equation*}
when no confusion is possible, we will omit the subscript \(\catA\) from
\(j\) and \({\mathcal A}\) from \(\ud\mu\), i.e., we will simply
write
\(\displaystyle{\int_{\mathcal A}j\ud\mu}\,\)
for
\(\displaystyle{\int_{\mathcal A}j^{}_\catA\ud\mu^{}_{\mathcal A}}\).
\par
The integral of \(j_\catA\) over \(\mathcal A\) is called the 
\emph{partition function} of
the groupoid \(\catA\).
Note that for any integrable morphism \(\varphi\colon\catA\to V\)
\begin{equation*}
\int_{\mathcal A}\varphi\,\ud \mu^{}_{\mathcal
A}=\varphi\left(\int_{\mathcal A}j_\catA\ud
\mu^{}_{\mathcal A}\right)
\end{equation*}
It is immediate to check that, for any two discrete
groupoids with finite homs and graded objects,
\begin{equation*}
\int_{\mathcal A{}^{\coprod}\mathcal B}j_{\catA{}^{\coprod}\catB}\ud
\mu^{}_{\mathcal A{}^{\coprod}\mathcal B}=\left(\int_{\mathcal
A}j_\catA\ud
\mu^{}_{\mathcal A}\right)\oplus\left(\int_{\mathcal B}j_\catB\ud
\mu^{}_{\mathcal B}\right)
\end{equation*}
and
\begin{equation*}
\int_{\mathcal A\times\mathcal B}j_{\catA\times\catB}\ud 
\mu^{}_{\mathcal A\times\mathcal
B}=\left(\int_{\mathcal A}j_\catA\ud \mu^{}_{\mathcal 
A}\right)\otimes\left(\int_{\mathcal
B}j_\catB\ud
\mu^{}_{\mathcal B}\right)\,.
\end{equation*}
If \(\varphi^{}_\catA\colon\catA\to V\) and \(\varphi^{}_\catB\colon\catB\to
W\) are two integrable morphisms, then also \(\varphi^{}_\catA\oplus
\varphi^{}_\catB\)  is integrable, and
\begin{align*}
\int_{\mathcal A{}^{\coprod}\mathcal B}(\varphi^{}_\catA\oplus
\varphi^{}_\catB)\ud
\mu^{}_{\mathcal A{}^{\coprod}\mathcal B}=\left(\int_{\mathcal
A}\varphi^{}_\catA\ud \mu^{}_{\mathcal
A}\right)\oplus\left(\int_{\mathcal B}\varphi^{}_\catB\ud
\mu^{}_{\mathcal B}\right)\,.
\end{align*}
Moreover, if \(\varphi^{}_\catA\) and \(\varphi^{}_\catB\) are graded,
then
\(\varphi^{}_\catA\otimes \varphi^{}_\catB\) is graded (therefore
integrable) and
\begin{align*}
\int_{\mathcal A\times\mathcal B}(\varphi^{}_\catA\otimes
\varphi^{}_\catB)\ud
\mu^{}_{\mathcal A\times\mathcal B}=\left(\int_{\mathcal
A}\varphi^{}_\catA\ud
\mu^{}_{\mathcal A}\right)\otimes\left(\int_{\mathcal
B}\varphi^{}_\catB\ud
\mu^{}_{\mathcal B}\right)
\,.
\end{align*}

A pull-back formula holds for integration over discrete groupoids.
\begin{lemma}
Let \(\pi\colon\catA\to \catB\) be a finite covering of discrete
groupoids with finite homs.
Then
\begin{equation*}
\int_{\mathcal A}\pi^*j_{\catB}\ud\mu^{}_{\mathcal A}=\int_{\mathcal
B}(\deg\pi\cdot j_{\catB})\ud\mu^{}_{\mathcal B}
\end{equation*}
In particular, for a homogeneous covering,
\begin{equation*}
\int_{\mathcal A}\pi^*j_{\catB}\ud\mu^{}_{\mathcal
A}=(\deg\pi)\int_{\mathcal B}j_{\catB}\ud\mu^{}_{\mathcal B}
\end{equation*}
\end{lemma}
\begin{proof}
By formulas \prettyref{eq:cardinality} and \prettyref{eq:degree-cov}, we have
\begin{align*}
\int_{\mathcal A}\pi^*j_{\catB}\ud\mu^{}_{\mathcal A}
&=\sum_{[x]\in\mathcal
A}\frac{[\pi(x)]}{\card{\Aut x}}
=\sum_{[y]\in\mathcal
B}\left(\sum_{[x]\in\pi^{-1}([y])
}\frac{[\pi(x)]}{\card{\Aut x}}\right)\\
&=\sum_{[y]\in\mathcal
B}\left(\sum_{[x]\in\pi^{-1}([y])
}\frac{1}{\card{\Aut x}}\right)[y]
=\sum_{[y]\in\mathcal
B}\left(\sum_{[x]\in\pi^{-1}([y])
}\card{[x]_y}\right)\frac{[y]}{\card{\Aut y}}\\
&=\sum_{[y]\in\mathcal
B}\deg\pi(y)\frac{[y]}{\card{\Aut y}}
\end{align*}
\end{proof}
As an immediate consequence we get
\begin{prop}[The pull-back formula]\label{prop:pull-back}
If \(\pi\colon\catA\to\catB\) is a homogeneous covering of discrete
groupoids with finite homs, then, for any integrable morphism
\(\varphi\colon\catB\to
V\) one has
\begin{equation*}
\boxed{
\int_{\mathcal A}\pi^*\varphi\,\ud\mu^{}_{\mathcal
A}=(\deg\pi)\int_{\mathcal B}\varphi\,\ud\mu^{}_{\mathcal B}}
\end{equation*}
\end{prop}
We end this section by introducing the concept of push-forward
of an integrable morphism and prove the Fubini theorem and the
push-pull formula in
this context.
\begin{dfn}
Let \(\pi\colon\catA\to\catB\) be a finite covering of discrete
groupoids with finite homs, and let \(\varphi\colon \catA\to V\) be an
integrable morphism. The \emph{push-forward} of \(\varphi\) is the morphism
\(\pi_*\varphi\colon\catB\to V\) defined by
\begin{equation*}
(\pi_*\varphi) (y)=\sum_{x\in {\rm Ob}\left((\pi^{-1}({\bf y})\right)} \varphi(x)
\end{equation*}
\end{dfn}

Note that the elements in the sum on the right-hand side are not weighted
by the factors \(1/\card{\Aut x}\), that is, the push-forward morphism
\(\pi_*\varphi\) is \emph{not} an integral over the fibre. This apparently
unnatural choice has two main motivations. The first is that the
combinatorial relation between a coordinate-free Feynman diagrams
expression and its reformulation in terms of a fixed system of coordinates
is a push-forward according to the above definition. The second is that we
want Fubini's theorem to hold.

\begin{prop}[Fubini's theorem]\label{prop:fubini} Let
\(\pi\colon\catA\to\catB\) be a finite covering of discrete
groupoids with finite homs, and let
\(\varphi\colon
\catA\to V\) be an integrable morphism. Then
\begin{equation}\label{eq:fubini}
\boxed{\int_{\mathcal B}\pi_*\varphi\,\ud\mu_{\mathcal B}=\int_{\mathcal
A}\varphi\,\ud\mu_{\mathcal A}}
\end{equation}
\end{prop}
\begin{proof}
 By definition of push-forward,
\begin{align*}
\int_{\mathcal B}\pi_*\varphi\,\ud\mu_{\mathcal B}&=
\sum_{[y]\in{\mathcal B}}\frac{1}{\card{\Aut y}}\left(\sum_{x\in {\rm
Ob}\left((\pi^{-1}({\bf y})\right)}
\varphi(x)\right)\\
&=
\sum_{[y]\in{\mathcal B}}\left(\sum_{[x]\in\pi^{-1}([y])}
\frac{\card{[x]_y}}{\card{\Aut
y}}\varphi(x)
\right).
\end{align*}
We now use equation \prettyref{eq:cardinality} to rewrite the right-hand
term as
\begin{align*}
\sum_{[y]\in{\mathcal B}}\left(\sum_{[x]\in\pi^{-1}([y])}
\frac{\varphi(x)}{\card{\Aut
x}}
\right)=\int_{\mathcal
A}\varphi\,\ud\mu_{\mathcal A}\,.
\end{align*}
\end{proof}
As a consequence, we find
\begin{prop}[The push-pull formula] Let \(\pi\colon\catA\to\catB\) be a
finite homogeneous
covering of discrete groupoids with finite homs, and let
\(\varphi\colon
\catA\to V\) be an integrable morphism. Then
\begin{equation}\label{eq:push-pull}
\boxed{\int_{\mathcal A}\pi^*\pi_*\varphi\,\ud\mu_{\mathcal
A}=(\deg \pi)\int_{\mathcal A}\varphi\,\ud\mu_{\mathcal A}}
\end{equation}

\end{prop}

\section{Symmetric powers of groupoids}
One of the basic equations of the theory of Feynman diagrams states
that exponentiation changes a sum over \emph{connected} diagrams
into a sum over all diagrams. The proof of this relation
is clearer in the general context of integration
over groupoids, so let us consider an arbitrary discrete
groupoid with finite homs
\(\catA\), and let
\({\sf S}_n\) be the simply connected groupoid having the
symmetric group \(\Perm{n}\)  as set of objects.  The \(n\)-th symmetric 
product of
a groupoid \(\catA\) with itself is the quotient groupoid
\(\text{Sym}_n(\catA)=\catA^{\times n}/\Perm{n}\).

More explicitly, the groupoid \(\text{Sym}_n(\catA)\) has the same objects
as
\(\catA^{\times n}\), but has additional morphisms given by the action 
of the permutation
group. For instance, for any object \(x\) of \(\catA\),
\begin{equation*}
\Aut_{\text{Sym}_n(\catA)}(\underbrace{x,x,\dots,x}_{n\text{
times}})=\Perm{n}\rtimes_\rho\prod_{i=1}^n\Aut{x}\,,
\end{equation*}
where the semi-direct product structure is induced by the natural map
\begin{equation*}
\rho\colon\Perm{n}\to \Aut\left(\prod_{i=1}^n\Aut{x}\right)\,.
\end{equation*}
The coarse space of \(\text{Sym}_n(\catA)\) is clearly 
\(\text{Sym}_n(\mathcal A)\). Note
that there is a natural immersion \(\iota\colon\catA^{\times
n}\hookrightarrow\text{Sym}_n(\catA)\). We set
\(\catA^\times\!:=\coprod_{n=0}^\infty \catA^{\times n}\) and
\(\text{\rm Sym}(\catA)\!:=\coprod_{n=0}^{\infty}\text{\rm
Sym}_n(\catA)\). Both group\-oids
\(\catA^\times\) and \(\text{\rm
Sym}(\sf A)\) have a symmetric monoidal category structure given by juxtaposition of
finite
sequences of objects of \(\catA\).

The map
\begin{align*}
\pi\colon \catA^{\times n}\times{\sf S}_n &\to \text{\rm Sym}_n(\catA)\\
(x_1,\dots,x_n,\sigma)&\mapsto (x_{\sigma(1)},\dots,x_{\sigma(n)})
\end{align*}
is a groupoid covering. For any \((x_1,\dots,x_n)\) in \(\text{\rm 
Sym}_n(\catA)\), the
objects of the fibre 
are the \((n+1)\)-ples
\((x_{\sigma(1)},\dots,x_{\sigma(n)},\sigma^{-1})\), which bijectively correspond
to the elements $\sigma$ in the symmetric group \(\Perm{n}\), so that 
\(\deg\pi=n!\). In
particular, it follows that
\begin{align}\label{eq:prod-vs-sym}
\notag\int_{\text{\rm Sym}_n(\mathcal A)} j_{\text{\rm Sym}_n(\catA)} \ud
\mu^{}_{\text{\rm Sym}_n(\mathcal A)}&=\frac{1}{n!}
\int_{{\mathcal A}^{\times n}\times {\mathcal S}_n}\pi^* j_{\text{\rm 
Sym}_n(\catA)} \ud
\mu^{}_{{\mathcal A}^{\times n}\times {\mathcal S}_n}\\
\notag&=\frac{1}{n!}\pi\left(\left(\int_{{\mathcal A}^{\times 
n}}j^{}_{\catA^{\times n}}\ud
\mu^{}_{{\mathcal A}^{\times n}}\right)\otimes\left(\int_{{\mathcal 
S}_n}j^{}_{{\sf S}_n}\ud
\mu^{}_{{\mathcal S}_n}\right)\right)\\
\notag&=\frac{1}{n!}\pi\left(\left(\int_{{\mathcal A}^{\times 
n}}j^{}_{\catA^{\times n}}\ud
\mu^{}_{{\mathcal A}^{\times n}}\right)\otimes[e]\right)\\
&=\frac{1}{n!}\iota\left(\int_{{\mathcal A}^{\times 
n}}j^{}_{\catA^{\times n}}\ud
\mu^{}_{{\mathcal A}^{\times n}}\right)\,,
\end{align}
where $e$ is the unit element in $\Perm{n}$ and the \(\iota\) in the last equation
is the immersion of
\(\catA^{\times n}\) in \(\text{\rm Sym}_n(\catA)\).

Let now \(R\) be a commutative topological \(\setC\)-algebra with unit,
 and let \(\varphi\colon\catA\to R\) be any morphism.
The product in \(R\) can be
seen as a tensor product on the trivial groupoid having the elements
of \(R\) as objects. Since
\(\catA^{\times}\) is freely generated by \(\catA\) as a symmetric monoidal
category,  \(\varphi\)
 uniquely extends to a tensor functor
\(\varphi\colon\catA^\times\to R\). Explicitly,
\begin{equation*}
\varphi(x_1,\dots,x_n)=\varphi(x_1)\cdots \varphi(x_n)
\end{equation*}
Due to the commutativity of \(R\), the restriction of \(\varphi\) to
\(\catA^{\times n}\) is
\(\Perm{n}\)-invariant, i.e.,  \(\varphi\) is a tensor
functor
\(\text{\rm Sym}({\sf A})\to R\). If \(\varphi\colon\catA\to R\) is
graded, then also \(\varphi\colon \text{\rm Sym}({\sf A})\to R\) is.
 Therefore, the diagram
\begin{equation*}
{\xymatrix{\displaystyle{
\catA^{\times
n}}&\displaystyle{
R}\\
\displaystyle{\text{\rm Sym}_n(\catA)}&
\ar "1,1";"1,2"^{\varphi}
\ar "2,1";"1,2"_{\varphi}
\ar "1,1";"2,1"_{\iota}
}}
\end{equation*}
commutes, and \prettyref{eq:prod-vs-sym} gives
\begin{equation}\label{eq:pre-part-vs-free}
\int_{\text{\rm Sym}_n(\mathcal A)} \varphi\, \ud
\mu^{}_{\text{\rm Sym}_n(\mathcal A)}=
\frac{1}{n!}\int_{{\mathcal A}^{\times n}}\varphi\,\ud
\mu^{}_{{\mathcal A}^{\times n}}\,.
\end{equation}

 Assume that the series
\(\sum_{n=0}^\infty r^n/n!\) is convergent in \(R\) and write
\begin{equation*}
\exp(r)=\sum_{n=0}^\infty\frac{r^n}{n!}\,.
\end{equation*}
With these notations, we have
\begin{prop}\label{prop:abstract-F-vs-Z}
 Let \(\catA\) be a discrete groupoid with finite homs having no objects
of degree zero,
and let \(R\) be any complete graded
\(\setC\)-algebra with unit. Then, for any graded functor
\(\varphi\colon\catA\to R\), the following
identity holds:
    \begin{equation}
      \label{eq:part-vs-free}
     \exp \left\{\int_{\mathcal A}\varphi\,\ud\mu^{}_{\mathcal
A}\right\}
      = \int_{\text{\rm Sym}(\mathcal A)} \varphi\, \ud
\mu^{}_{\text{\rm Sym}(\mathcal A)}\,.
    \end{equation}
\end{prop}
\begin{proof}
    By definition of exponential and using 
\prettyref{eq:pre-part-vs-free} we find
    \begin{align*}
\sum_{n=0}^\infty \oneof{n!} \left(\varphi\int_{\mathcal
A}j_\catA^{}\ud\mu^{}_{\mathcal A}\right)^{n}
&= \sum_{n=0}^\infty \oneof{n!} \varphi\left(\int_{\mathcal
A}j_\catA^{}\ud\mu^{}_{\mathcal A}\right)\tp{n}\\
&= \sum_{n=0}^\infty \oneof{n!} \varphi\left(\int_{\mathcal A^{\times
n}}j_{\catA^{\times n}}^{}\ud\mu^{}_{\mathcal A^{\times n}}\right)\\
  &= \sum_{n=0}^\infty \oneof{n!} \int_{\mathcal A^{\times
n}}\varphi\,\ud\mu^{}_{\mathcal A^{\times n}}\\  
&=\sum_{n=0}^\infty\int_{\text{\rm Sym}_n(\mathcal A)} \varphi\, \ud
\mu^{}_{\text{\rm Sym}_n(\mathcal A)}\,.
    \end{align*}
   \end{proof}

\section{Feynman diagrams}\label{sec:feynman-diagrams}
 A Feynman diagram is, basically, a graph with a distinguished subset
of 1-valent vertices (the endpoints) and some additional
structure on the other vertices (the internal vertices). Moreover Feynman
diagrams can be given a splitting of the endpoints into \emph{inputs} and
\emph{outputs}.
The additional structure on the internal vertices are a colour
decorating the vertex and a combinatorial datum adding some ``rigidity'',
i.e., reducing the automorphism group of the germs of edges stemming from
a vertex (the legs of the vertex) to a proper subgroup of the symmetric
group. In full generality,
this combinatorial datum could be any structure on the finite set of
the legs of a given vertex; on the other hand, the essence of
Feynman diagrams is their being a powerful graphical tool for computing
asymptotic expansions of Gaussian integrals. Therefore, we will consider
only those combinatorial data that the legs of a vertex ``naturally
inherits by drawing a graph on a sheet of paper''. 
More specifically, we will consider only the following three kind of
vertices: \emph{coupon vertices}, whose legs are split 
in the two subsets of legs stemming on the top and on the bottom of
the vertex\footnote{The reader familiar with Joyal-Street's tensor
calculus \cite{joyal-street;1991} or Reshetikhin-Turaev's graphical
calculus \cite{reshetikhin-turaev;ribbon-graphs,bakalov-kirillov} will 
find this
splitting familiar.}, \emph{cyclic vertices} whose legs inherit a 
cyclic order by the orientation of 
the sheet of paper, and \emph{symmetric vertices} whose legs inherit no 
additional structure, so that the automorphism group of a symmeric
vertex is the full symmetric group on its legs. 
Finally, any vertex has
an additional decorating \emph{colour}; the sets of colours we can put on
coupon, cyclic and symmetric vertices will be denoted
by the symbols {\sf Co}, {\sf Cy}, and {\sf Sy} respectively. In general, the
colouring sets can vary with the valence, for instance we may have just one 
colour for 3-valent vertices, and 37 (or even infinite) colours for 4-valent
vertices.  
\par
In the theory of Feynman diagrams one often deals with huge
families of colours decorating the vertices,\footnote{Actually, one 
encounters also
colours decorating the edges. Here we preferred not to work in full 
generality here
to make the exposition clearer. However, a short digression on Feynman 
diagrams with
coloured edges will be done later, relating the colours on the edges 
with the
elements of a chosen basis in a given vector space.} but only few of
them generally occur.  To formalize this feature, each set of colours is 
split
into the disjoint subsets of
\emph{ordinary} and
\emph{special} colours. Moreover, given a Feynman diagram, it is 
important to
have the possibility to distinguish one of its ordinary vertices, and
to look at it as a special one. To do this, an identification
of the set of ordinary colours with a subset of the special ones is given;
this allows us to talk of the special colour corresponding to a given
ordinary colour. Colours decorating the vertices will be denoted
by Greek letters; the special colour corresponding to the ordinary colour
$\alpha$ will be denoted by $\bmath{\alpha}$.
We will further assume that, for any valence \(d\),
the set of ordinary colours for the \(d\)-valent vertices is finite. In many
applications there is just one special and one ordinary colour for each valence,
so that the datum of the coloring is actually redundant and is reduced to
 the label ``ordinary'' or ``special''.
\par
 A  formal
definition of Feynman diagram is the following.
\begin{dfn}
Let \(n\) be a natural number. A Feynman diagram with $n$ legs
is the following set of data:
\begin{enumerate}
\item a 1-dimensional CW-complex \(\Gamma\);
\item \(n\) distinguished 1-valent vertices of \(\Gamma\), called the
\emph{endpoints} of \(\Gamma\); the vertices of $\Gamma$ which are not
endpoints will be called \emph{internal vertices}. The germs of the edges
stemming
from the endpoints are called \emph{legs} of the diagram and are denoted
by \(\text{Legs}(\Gamma)\); any edge which does not end in an
endpoint will be called an \emph{internal edge} of the diagram;
\item
a map
\begin{equation*}
\text{Internal vertices}(\Gamma)\to
\text{\sf
Co}\cup
\text{\sf Cy}\cup\text{\sf Sy}\,
\end{equation*}
called \emph{decoration};

\item elements in the pre-image of $\sf Co$ are the \emph{coupon} 
vertices; for
any \emph{coupon} vertex
\(v\) of
\(\Gamma\), a splitting of the set
\(\Leg(v)\) --- i.e., germs of edges stemming from \(v\) ---  into two
totally
ordered subsets denoted \(\In(v)\) and
\(\Out(v)\) respectively. A coupon vertex decorated by the colour \(\alpha\)
is depicted by
\[\xy*!LC\xybox{
        (0,1)*+[F]{\alpha};%
        (-.8,0)**\dir{-},(-0.4,0)**\dir{-},%
        (0.2,0.1)*+{\ldots},(.8,0)**\dir{-},%
        (-.8,2)**\dir{-},(-0.4,2)**\dir{-},%
        (0.2,1.8)*+{\ldots},(.8,2)**\dir{-}%
        }\endxy\,,
\]
with inputs on the lower side and outputs on the upper side, the total
order being that induced by the horizontal coordinate in the plane;
\item elements in the pre-image of $\sf Cy$ are the \emph{cyclic} 
vertices; for
any \emph{cyclic} vertex \(v\) of \(\Gamma\), a cyclic order on the set
\(\Leg(v)\). A cyclic vertex decorated by the colour \(\beta\) is
depicted by
\[\xy*!LC\xybox{
\rgvertex[\scriptstyle{\beta}]{7}%
\loose1\loose2\loose3\loose4\missing5\loose6%
\loose7
        }\endxy\,,
\]
the cyclic order being that inherited by the standard
counterclockwise orientation on the plane;

\item elements in the pre-image of $\sf Sy$ are the \emph{symmetric} 
vertices; no
additional combinatorial constraints are imposed on a symmetric vertex; a
symmetric vertex decorated by the colour
\(\gamma\) is depicted by
\[\xy*!LC\xybox{
        \blackvertex{7}\loose1\loose2\loose3\loose4\missing5\loose6%
\loose7
        }
,(1.1,0.3)*{\scriptstyle{\gamma}}\endxy
\]
\end{enumerate}
The set of colours are split into the two disjoint subsets of \emph{ordinary} and
\emph{special} colours. A vertex will be called ordinary or special
depending on the colour decorating it. Moreover, an injective map from ordinary to
special colours is given and, for any valence $d$, the set of colours
decorating $d$-valent ordinary vertices is finite.
Isomorphisms  between
Feynman diagrams are homotopy classes of cellular isomorphisms of the
underlying CW-complexes, respecting the additional structures on the
vertices.
\end{dfn}
 It is clear from this definition that Feynman diagrams are deeply related
with the graphical formalism of \cite{joyal-street;1991} and
\cite{reshetikhin-turaev;ribbon-graphs}.
If the only morphisms we allow between Feynman diagrams are
the isomorphisms, then Feynman diagrams form a groupoid. The groupoid of
Feynman
diagrams with $n$ legs will be denoted by the symbol
\(\catF(n)\).
It is immediate to check that the groupoids
\(\catF(n)\) are discrete groupoids with finite homs.
Moreover, the objects of the groupoids \(\catF(n)\) are
\emph{graded} by the sum of the valencies of the ordinary vertices.

We now add a splitting of the endpoints into ``inputs'' and ``outputs''

\begin{dfn}
Let \(n\) and \(m\) be two natural numbers. A Feynman diagram of type
\((n,m)\) is the following set of data:
\begin{enumerate}
\item a Feynman diagram $\Gamma$ with $(m+n)$ legs;
\item a splitting of \(\text{Endpoints}(\Gamma)\) into two
subsets \(\In(\Gamma)\) and \(\Out(\Gamma)\), of cardinality \(n\) and
\(m\) respectively. The set \(\In(\Gamma)\) is called the set of the
\emph{inputs} of \(\Gamma\) and the set \(\Out(\Gamma)\) is called the set
of the \emph{outputs} of \(\Gamma\);
\item two bijections
\begin{align*}
\iota^{}_\Gamma\colon\{1_{\text{in}},
\dots,n_{\text{in}}\}&\to\In(\Gamma)\\
\omega^{}_\Gamma\colon\{1_{\text{out}},
\dots,m_{\text{out}}\}&\to\Out(\Gamma)
\end{align*}
inducing total orders on \(\In(\Gamma)\) and \(\Out(\Gamma)\).
\end{enumerate}
\end{dfn}
Isomorphisms between Feynman diagrams of type $(n,m)$ are isomorphisms of
the underlying Feynman diagrams with $(m+n)$ legs respecting the
additional structure. The groupoid of Feynman diagrams
of type
\((n,m)\),  will be denoted by the symbol
\(\catF(n,m)\). There is a natural ``forget the numbering'' functor
\begin{equation*}
\catF(n,m)\to\catF(m+n).
\end{equation*}
If $n=m=0$, this is an isomorphism
\begin{equation*}
\catF(0,0)\to\catF(0).
\end{equation*}
\par
We now discuss some 2-categorical features of Feynman diagrams with inputs
and outputs. 

\par
If $\Gamma$ is a Feynman diagram of type $(n,m)$, the natural numbers
\(n\) and \(m\) are called the \emph{source} and the
\emph{target} of the Feynman diagram \(\Gamma\) and are denoted by
the symbols \(\Src(\Gamma)\) and \(\Tgt(\Gamma)\).
If $\Gamma$ and $\Phi$ are two Feynman diagrams
such that
$\Src(\Gamma) = \Tgt(\Phi)$, then we
can
form a new Feynman diagram $\Gamma\circ\Phi$ by identifying the point
\(\nu_{\text{in}}\) of \(\Gamma\) with the point \(\nu_{\text{out}}\)
of \(\Phi\) for any \(\nu\) between \(1\) and \(\Src{\Gamma}=
\Tgt{\Phi}\); we have \(\In(\Gamma\circ\Phi):=\In(\Phi)\) and
\(\Out(\Gamma\circ\Phi):=\Out(\Gamma)\). An example is
\begin{equation*}
{\xy*!LC\xybox{(0,.5)*{\bullet};%
        (0,0)**\dir{-},(-.6,1.1)**\dir{-},(.6,1.1)**\dir{-}%
        ,(-.3,0.4)*{\scriptstyle{\alpha}}%
,(0,-.2)*{1_{\text{\rm in}}},(-.7,1.3)*{1_{\text{\rm out}}}%
        ,(.7,1.3)*{2_{\text{\rm out}}}%
        }\endxy}\quad\circ\quad
{\xy*!LC\xybox{
        (0,.5)*+[F]{\beta};%
        (-.6,-.35)**\dir{-},%
        (.6,-.35)**\dir{-},%
        (0,1.5)*{\,\,\gamma\,\,}*\cir{}**\crv{(-0.2,1) & (-1,1.2)},%
        (0,1.5)*{\,\,\gamma\,\,}*\cir{}**\crv{(.2,1) & (1,1.2)},
        (0,1.5)*{\,\,\gamma\,\,}*\cir{};(0,2.1)**\dir{-}%
        ,(0,2.3)*{1_{\text{\rm out}}},(-.6,-.55)*{1_{\text{\rm in}}}%
        ,(.7,-.55)*{2_{\text{\rm in}}}%
        }\endxy}\quad=\quad
{\xy*!LC\xybox{
        (0,.5)*+[F]{\beta};%
        (-.6,-.35)**\dir{-},%
        (.6,-.35)**\dir{-},%
        (0,1.5)*{\,\,\gamma\,\,}*\cir{}**\crv{(-0.2,1) & (-1,1.2)},%
        (0,1.5)*{\,\,\gamma\,\,}*\cir{}**\crv{(.2,1) & (1,1.2)},
        (0,1.5)*{\,\,\gamma\,\,}*\cir{};(0,2.4)*{\bullet}**\dir{-}%
        ,(0,2.4);(-.6,3)**\dir{-},(.6,3)**\dir{-}
        ,(-.3,2.3)*{\scriptstyle{\alpha}}
        ,(-.6,-.55)*{1_{\text{\rm in}}}%
        ,(.7,-.55)*{2_{\text{\rm in}}}%
        ,(-.7,3.2)*{1_{\text{\rm out}}}%
        ,(.7,3.2)*{2_{\text{\rm out}}}%
        }\endxy}
\end{equation*}
In this way we have defined a \emph{composition}
\begin{equation*}
\catF(n,m)\times \catF(k,n)\to \catF(k,m)\,.
\end{equation*}
The right way to look at this composition of diagrams is as a category
structure on the set of natural numbers where Feynman diagrams of type
\((n,m)\) are the morphisms between \(n\) and \(m\) (see
\cite{baez-dolan;finite-sets}). Moreover, since
\(\catF(n,m)\) is a groupoid, we have defined a category whose Hom-spaces
are groupoids, i.e., if we denote by the symbol \(\catN\) the trivial
groupoid having the natural numbers as objects, then
\begin{equation*}
\catF\colon \catN\times\catN\to {\sf Grpd}\,,
\end{equation*}
 is an \emph{enriched category} in the
sense of Eilenberg-Kelly
(\cite{eilenberg-kelly,kelly;enriched-category}). In particular,
$\catF$ is a 2-category,\footnote{The reader familiar with higher dimensional
category theory will immediately notice that composition in
$\catF$  is not strictly associative and also identity arrows are not strict, being
so only up to a natural isomorphism. That is, $\catF$ is a bicategory (or weak
2-category) rather than a 2-category. On the other hand, the associativity and
unital constrains of the bicategory $\catF$ are quite evident, so we preferred not
to specify them aiming to a sort of compromise between the complete rigour of
higher  dimensional category theory and an exposition enjoyable by the non-expert
reader.} with the natural numbers as objects. The identity morphism
\(j_n\colon n\to n\) is clearly given by
\begin{equation*}
j_n=
{\xy*!LC\xybox{
        (-1,0);(-1,1)**\dir{-}
        ,(-.2,0);(-.2,1)**\dir{-}
        ,(1.5,0);(1.5,1)**\dir{-}
        ,(-1,-.2)*{1_{\text{\rm in}}}%
        ,(-1,1.2)*{1_{\text{\rm out}}}%
        ,(-.2,-.2)*{2_{\text{\rm in}}}%
        ,(-.2,1.2)*{2_{\text{\rm out}}}%
        ,(1.5,-.2)*{n_{\text{\rm in}}}%
        ,(1.5,1.2)*{n_{\text{\rm out}}}%
        ,(.6,.5)*{\cdots}
        }\endxy}
\end{equation*}

 Given any two Feynman diagrams
\(\Gamma\) and \(\Phi\), we can make a new Feynman diagram out of them
 by taking their disjoint union. Some care has to be taken in
defining the new numberings on the inputs and the outputs: if we
denote the new Feynman diagram by the symbol  
\(\Gamma\otimes \Phi\), then
\begin{enumerate}
\item \(\text{Endpoints}(\Gamma\otimes
\Phi):=\text{Endpoints}(\Gamma)\cup \text{Endpoints}(\Phi)\);
\item \(\In(\Gamma\otimes \Phi):=\In(\Gamma)\cup
\In(\Phi)\);
\item \(\Out(\Gamma\otimes \Phi):=\Out(\Gamma)\cup
\Out(\Phi)\);
\item
\begin{equation*}
\iota_{\Gamma\otimes\Phi}(\nu)=\left\{
\begin{matrix}
\iota_\Gamma(\nu)&\text{if}&1\leq\nu\leq \Src{\Gamma}\\
\iota_\Phi(\nu-\Src(\Gamma))&\text{if}&\Src(\Gamma)+1\leq\nu\leq
\Src{\Gamma}+\Src{\Phi}
\end{matrix}\right.
\end{equation*}
\item
\begin{equation*}
\omega_{\Gamma\otimes\Phi}(\nu)=\left\{
\begin{matrix}
\omega_\Gamma(\nu)&\text{if}&1\leq\nu\leq \Tgt{\Gamma}\\
\omega_\Phi(\nu-\Tgt(\Gamma))&\text{if}&\Tgt(\Gamma)+1\leq\nu\leq
\Tgt{\Gamma}+\Tgt{\Phi}
\end{matrix}\right.
\end{equation*}
\end{enumerate}
For instance,
\begin{equation*}
{\xy*!LC\xybox{(0,.5)*{\bullet};%
        (0,0)**\dir{-},(-.6,1.1)**\dir{-},(.6,1.1)**\dir{-}%
        ,(-.3,0.4)*{\scriptstyle{\alpha}}%
,(0,-.2)*{1_{\text{\rm in}}},(-.7,1.3)*{1_{\text{\rm out}}}%
        ,(.7,1.3)*{2_{\text{\rm out}}}%
        }\endxy}\quad\otimes\quad
{\xy*!LC\xybox{
        (0,.5)*+[F]{\beta};%
        (-.6,-.35)**\dir{-},%
        (.6,-.35)**\dir{-},%
        (0,1.5)*{\,\,\gamma\,\,}*\cir{}**\crv{(-0.2,1) & (-1,1.2)},%
        (0,1.5)*{\,\,\gamma\,\,}*\cir{}**\crv{(.2,1) & (1,1.2)},
        (0,1.5)*{\,\,\gamma\,\,}*\cir{};(0,2.1)**\dir{-}%
        ,(0,2.3)*{1_{\text{\rm out}}},(-.6,-.55)*{1_{\text{\rm in}}}%
        ,(.7,-.55)*{2_{\text{\rm in}}}%
        }\endxy}\quad=\quad
{\xy*!LC\xybox{(0,.5)*{\bullet};%
        (0,0)**\dir{-},(-.6,1.1)**\dir{-},(.6,1.1)**\dir{-}%
        ,(-.3,0.4)*{\scriptstyle{\alpha}}%
,(0,-.2)*{1_{\text{\rm in}}},(-.7,1.3)*{1_{\text{\rm out}}}%
        ,(.7,1.3)*{2_{\text{\rm out}}}%
        }\endxy}
{\xy*!LC\xybox{
        (0,.5)*+[F]{\beta};%
        (-.6,-.35)**\dir{-},%
        (.6,-.35)**\dir{-},%
        (0,1.5)*{\,\,\gamma\,\,}*\cir{}**\crv{(-0.2,1) & (-1,1.2)},%
        (0,1.5)*{\,\,\gamma\,\,}*\cir{}**\crv{(.2,1) & (1,1.2)},
        (0,1.5)*{\,\,\gamma\,\,}*\cir{};(0,2.1)**\dir{-}%
        ,(0,2.3)*{3_{\text{\rm out}}},(-.6,-.55)*{2_{\text{\rm in}}}%
        ,(.7,-.55)*{3_{\text{\rm in}}}%
        }\endxy}
\end{equation*}
The choice of the symbol \(\Gamma\otimes \Phi\) to denote this graph is
not
an accident: it actually is a tensor product of morphisms in the enriched
category \(\catF\), which is therefore a monoidal 2-category (see, for
instance, \cite{kapranov-voevodsky,day-street}).

Moreover the Feynman diagram
\begin{equation*}
\sigma_{m,n}^{}={\xy*!LC\xybox{
        (-1.3,0);(-1.3,1)**\dir{-}
        ,(0,0);(0,1)**\dir{-}
        ,(2,0);(2,1)**\dir{-}
         ,(3.3,0);(3.3,1)**\dir{-}
        ,(4.6,0);(4.6,1)**\dir{-}
        ,(7,0);(7,1)**\dir{-}
        ,(1,.5)*{\cdots}
        ,(-1.3,-.2)*{1_{\text{\rm in}}}%
        ,(-1.3,1.2)*{{n+1}_{\text{\rm out}}}%
        ,(0,-.2)*{2_{\text{\rm in}}}%
        ,(0,1.2)*{{n+2}_{\text{\rm out}}}%
        ,(2,-.2)*{m_{\text{\rm in}}}%
        ,(2,1.2)*{{n+m}_{\text{\rm out}}}%
        ,(3.3,-.2)*{{m+1}_{\text{\rm in}}}%
        ,(3.3,1.2)*{1_{\text{\rm out}}}%
        ,(4.6,-.2)*{{m+2}_{\text{\rm in}}}%
        ,(4.6,1.2)*{2_{\text{\rm out}}}%
        ,(7,-.2)*{{m+n}_{\text{\rm in}}}%
        ,(7,1.2)*{{n}_{\text{\rm out}}}%
        ,(5.8,.5)*{\cdots}
        }\endxy}
\end{equation*}
is  a
\emph{braiding} operator between the tensor product of \(m\) with \(n\) and
the tensor product of \(n\) with \(m\). This is more evident if
one
draws $\sigma_{m,n}$ as follows:
\[
\sigma_{m,n}={\xy*!LC\xybox{(0,0);(-3,2)**\dir{-}
,(1.2,0);(-1.8,2)**\dir{-}
,(3,0);(0,2)**\dir{-}
,(-3,0);(-1.5,.75)**\dir{-}
,(-1,0);(-.5,.25)**\dir{-}
,(.75,1.85);(1,2)**\dir{-}
,(1.5,1.25);(3,2)**\dir{-}
,(.6,1)*{\cdots}
,(0.1,-.2)*{{m+1}_{\text{\rm in}}},(1.4,-.2)*{{m+2}_{\text{\rm in}}}%
,(-3,-.2)*{{1}_{\text{\rm in}}},(-1,-.2)*{{m}_{\text{\rm in}}}%
        ,(3.2,-.2)*{{m+n}_{\text{\rm in}}}%
,(-3,2.2)*{{1}_{\text{\rm out}}},(-1.8,2.2)*{{2}_{\text{\rm out}}}%
,(1.1,2.2)*{{n+1}_{\text{\rm out}}},(3.2,2.2)*{{n+m}_{\text{\rm out}}}%
        ,(-0.1,2.2)*{{n}_{\text{\rm out}}}%
,(-1.5,.35)*{\cdots},(1.5,1.75)*{\cdots}
        }\endxy}
\]
For any $k,m,n\in\setN$ and any $\Gamma\in \catF(k,m)$, there are evident
isomorphisms
\begin{align*}
\sigma_{m,n}\circ (\Gamma\otimes j_n)&\simeq (j_n\otimes \Gamma)\circ 
\sigma_{k,n}\\
\sigma_{n,m}\circ (j_n \otimes \Gamma)&\simeq (\Gamma\otimes j_n)\circ
\sigma_{n,k}\\
\end{align*}
\vskip -32pt
\begin{align*}
(\sigma_{k,m}\otimes j_n)\circ(j_k\otimes\sigma_{m,n})
&\simeq \sigma_{k+m,n}\\
(j_m\otimes\sigma_{k,n})\circ(\sigma_{k,m}\otimes
j_n)&\simeq
\sigma_{k,m+n}
\end{align*}
which satisfy the axioms of 2-braidings 
\cite{kapranov-voevodsky}.
Moreover, there is an evident isomorphism of Feynman diagrams
\[
\sigma_{n,m}\circ \sigma_{m,n}\simeq j_{m+n}\,,
\]
so that the braiding \(\sigma\) is symmetric. Therefore $\catF$ is a
symmetric monoidal 2-category having the natural numbers as objects; such
a structure can be
called 2-PROP.\footnote{The reader unfamiliar with the notion
of PROP may think to it just as a symmetric monoidal category having the natural
numbers as objects: this quite imprecise definition suffices to deal with the use
of PROPs done in this paper; see
\cite{adams;infinite-loop-spaces} for a formal definition of PROP.}
\par
As a remark, note that any braiding \(\sigma_{m,n}\) can be obtained
as a composition of
elementary braidings of the form
\begin{equation*}
j^{}_\mu\otimes\sigma^{}_{1,1}\otimes
j_\nu^{}\,,
\end{equation*}
with \(\mu+\nu=m+n-2\).
\par
Finally, if \(\Gamma\) and \(\Phi\) are two Feynman diagrams, then
\(\deg(\Gamma\otimes \Phi)=\deg(\Gamma)+\deg(\Phi)\) and, if the
composition \(\Gamma\circ\Phi\) is defined, \(\deg(\Gamma\circ
\Phi)=\deg(\Gamma)+\deg(\Phi)\), i.e., the 2-PROP \(\catF\) is a
graded 2-PROP.

\section{Feynman diagrams with distinguished
sub-diagrams}\label{sec:Fey-dia-dist}

Sums over graphs in the theory of Feynman diagrams are usually sums 
over Feynman
diagrams containing a distinguished sub-diagram. The aim of this section 
is to
formally introduce the concept of Feynman diagram with a distinguished 
sub-diagram,
and to discuss how this concept is related to coverings of groupoids.
\par
We have seen in \prettyref{sec:feynman-diagrams} that Feynman diagrams 
with inputs
and outputs can be seen as the 1-morphisms of a 2-category; in particular a
composition of Feynman diagrams with inputs
and outputs is defined. Also recall that there is a natural ``forget
the 
numbering''
morphism $\catF(m,n)\to\catF(m+n)$, which induces an isomorphism 
$\catF(0,0)\simeq
\catF(0)$. It is immediate to check that the projection $\catF(m,n)\to\catF(m+n)$ is
a homogeneous groupoid covering of degree $(m+n)!$: the objects in the fibre over a
Feynman diagram $\Gamma$ are the $(m+n)!$ Feynman diagrams of type $(m,n)$ which are
obtained by numbering the endpoints of $\Gamma$ in all possible ways, and the
isomorphisms in $\catF(m+n)$ uniquely lift to isomorphisms in $\catF(m,n)$.

\begin{dfn}
Let \(\Gamma\in\catF(n)\) be a Feynman diagram with \(n\) legs. We
say that \(\Gamma\) is a sub-diagram of the Feynman diagram 
\(\Psi\in\catF(m)\)
if there is an isomorphism of Feynman diagrams \(\hat\Psi\simeq
\hat\Gamma\circ\hat\Phi\) for a suitable pre-image
\(\hat\Gamma\) of \(\Gamma\) in \(\catF(n,0)\), a suitable pre-image
\(\hat\Psi\) of \(\Psi\) in \(\catF(m,0)\) and some Feynman diagram
\(\hat\Phi\in\catF(m,n)\).
\end{dfn}
\begin{dfn}
Let \(\Gamma\) be a Feynman diagram.
The groupoid \(\catF_{\Gamma}\) is the subgroupoid of \(\catF\) whose 
objects
 are the Feynman diagrams
\(\Psi\) such that
\begin{enumerate}
\item \(\Gamma\) is a distinguished sub-diagram of \(\Psi\);
\item all the vertices of \(\Psi\) outside \(\Gamma\) are ordinary;
\end{enumerate}
 Morphisms between objects of \(\catF_{\Gamma}\) are the morphisms in
\(\catF\) which map the distinguished sub-diagram \(\Gamma\) to itself.
\end{dfn}

Note that \(\catF_\emptyset\) is the groupoid of Feynman diagrams with only
ordinary vertices.
The objects of the groupoids
\(\catF_{\Gamma}\) are graded by the sum of the
 valencies of the ordinary vertices. Note that, for any $\Gamma$ and any fixed
degree
\(d\) there are only finitely many isomorphism classes of degree
\(d\) objects of \(\catF_{\Gamma}\).
\par

Observe that,
since any two pre-images of \(\Gamma\) in \(\catF(n,0)\) may only differ 
by the
action of an element of the symmetric group \(\Perm{n}\) on the
the inputs,  any Feynman diagram in \(\catF_{\Gamma}(0)\) can be
obtained as the composition of \emph{any}
pre-image of
\(\Gamma\) in \(\catF(n,0)\) with a suitable Feynman diagram in
\(\catF_{\emptyset}(0,n)\). Namely, the following proposition
holds.

\begin{prop}\label{prop:homo-cov-dia} Let $\Gamma$ be a Feynman diagram with $n$ legs and denote by
$\pi^{-1}(\bmath{\Gamma})$ the fibre over
$\bmath{\Gamma}$ in the projection $\pi: \catF(n,0)\to \catF(n)$. Then the
composition
$\circ\colon\pi^{-1}({\bmath \Gamma})\times
\catF_\emptyset(0,n)\to\catF_{\Gamma}(0)$ is a homogeneous groupoid covering of
degree
$n!$.
\end{prop}

\begin{proof}
As remarked above, the projection $\pi: \catF(n,0)\to \catF(n)$
 is a covering of degree $n!$ and the objects in
$\pi^{-1}(\bmath{\Gamma})$ are the same CW-complex underlying $\Gamma$, with the
additional datum of a numbering on endpoints.
The composition of diagrams induces a map
$\circ: \pi^{-1}(\bmath{\Gamma}) \times \catF_\emptyset(0,n) \to \catF_\Gamma(0)$,
which is a degree $n!$ covering, too. Indeed,
 let $\varphi: \Psi_1\to \Psi_2$ be a morphism in $\catF_\Gamma(0)$, and let
$(\Gamma_1,\Phi_1)$ be a pre-image of $\Psi_1$ in $\pi^{-1}(\bmath{\Gamma}) \times
\catF_\emptyset(0,n)$. The pre-image $(\Gamma_1,\Phi_1)$ is simply obtained by
``cutting'' $\Psi_1$ along $\Gamma$ and numbering the endpoints on the two
pieces in a compatible way.
Since $\varphi$ is a morphism in $\catF_\Gamma(0)$, it is a  Feynman diagrams
isomorphism which preserves $\Gamma$. So it will induce an automorphism of
$\Gamma$ and an isomorphism between $(\Psi_1\setminus\Gamma)$ and
$(\Psi_2\setminus \Gamma)$, where we are looking at these graphs as objects of
$\catF(n)$.    Since forgetting the numbering on endpoints is a covering,
these two $\catF(n)$-morphisms can be lifted (uniquely) to a
$\pi^{-1}(\bmath{\Gamma})$-morphism
$\Gamma_1 \to\Gamma_2$, and to a $\catF_\emptyset(0,n)$-morphism $\Phi_1 \to
\Phi_2$, i.e., the $\catF_\Gamma(0)$-morphism $\varphi: \Psi_1\to \Psi_2$ can
be lifted to a $(\pi^{-1}(\bmath{\Gamma}) \times \catF_\emptyset(0,n))$-morphism
$(\Gamma_1,\Phi_1)\to(\Gamma_2,\Phi_2)$. Therefore $\circ: \pi^{-1}(\bmath{\Gamma})
\times
\catF_\emptyset(0,n) \to \catF_\Gamma(0)$ is a fibration. Moreover the lifting
is unique, so $\circ$ is a covering.
\end{proof}

\everyxy={0,<2em,0em>:,(0,0.5),} 

We are now ready to give a formal treatment of the following basic principle in
the combinatorics of Feynman diagrams: if $\catF_1$,
$\catF_2$ and $\catF_3$ are three families of Feynman diagrams such that 
any diagram in
$\catF_1$ is the composition of a diagram in $\catF_2$ with a diagram in 
$\catF_3$, then
summing over $\catF_1$ is the same thing as summing over $\catF_2$ and 
$\catF_3$ and then
composing the two sums:
\begin{equation*}
\int_{\Fey_1}j\ud\mu=
\left(\int_{\Fey_2}j\ud\mu\right)\circ\left(\int_{\Fey_3}j\
\ud\mu\right).
\end{equation*}
For instance, if \(\Gamma\) is a Feynman diagram with \(n\) legs, then any
Feynman diagram in
$\catF_\Gamma(0)$ can be obtained by the composition of a preimage of
$\Gamma$ in $\catF(n,0)$ with some Feynman diagram in 
$\catF_\emptyset(0,n)$. This
implies
\begin{equation}\label{eq:first-gr-avg-two}
\int_{\Fey_{\Gamma}(0)}j\ud\mu=\left(\int_{\Gamma}j\ud\mu
\right)
\circ\left(\int_{\Fey_\emptyset(0,n)}j\ud\mu\right)
\end{equation}
The proof of \prettyref{eq:first-gr-avg-two} is almost immediate: from
\prettyref{prop:homo-cov-dia} and \prettyref{prop:pull-back} we know that
\begin{equation*}
\int_{\Fey_{\Gamma}(0)}j\ud\mu=\frac{1}{n!}
\left(\int_{\pi^{-1}(\Gamma)}j\ud\mu
\right)
\circ\left(\int_{\Fey_\emptyset(0,n)}j\ud\mu\right)
\end{equation*}
Due to the lack of an ordering on the legs of
\(\Gamma\), a composition with
\(\Gamma\) is not well-defined. On the other hand, any two
pre-images of \(\Gamma\) in \(\catF(n,0)\) may only differ by the
action of an element of the symmetric group \(\Perm{n}\). Therefore a
multiplication by \(\Gamma\) is well defined on the space of
\(\Perm{n}\)-invariant elements of \(\propF(0,n)\).\footnote{We are using the
notations from \prettyref{sec:int-over-grpd}: if $\catA$ is a groupoid, then
${\mathcal A}$ denotes the set of isomorphism classes of objects in $\catA$ and
${\mathscr A}$ denotes the free $\setC$-vector space generated by ${\mathcal A}$.}
Thus, we have an operator
\begin{equation}\label{eq:multiplication}
m^{}_{\Gamma}\colon \propF(0,n)^{\Perm{n}}\to
\propF(0)\,.
\end{equation}
In particular,
\begin{equation*}
m^{}_{\Gamma}\colon \propF_\emptyset(0,n)^{\Perm{n}}\to
\propF_{\Gamma}(0)\,.
\end{equation*}
Extending \prettyref{eq:multiplication} by linearity one obtains a
composition
\begin{equation}
\circ\colon \propF(n)\otimes\propF(0,n)^{\Perm{n}}\to
\propF(0)\,.
\end{equation}
The partition function
\begin{equation*}
\int_{\Fey_\emptyset(0,n)}j\ud\mu
\end{equation*}
is clearly \(\Perm{n}\)-invariant, so the composition on the right-hand side of
\prettyref{eq:first-gr-avg-two} 
is well defined. If $\hat\Gamma$ is any object in $\pi^{-1}({\bmath{\Gamma}})$, then
the operators $m_\Gamma$ and $m_{\hat{\Gamma}}$ on the space of
$\Perm{n}$-invariants of $\propF_\emptyset(0,n)$ do coincide, so that we have a
commu\-tative diagram
\begin{equation*}
{\xymatrix{\displaystyle{\pi^{-1}({\bmath{\Gamma}})}&\displaystyle{
\Hom(\propF_\emptyset(0,n)^{\Perm{n}},
\propF_\Gamma(0))}\\
\displaystyle{{\bmath{\Gamma}}}&
\ar "1,1";"1,2"^{m\phantom{mmmmmi}}
\ar "2,1";"1,2"_{m}
\ar "1,1";"2,1"_{\pi}
}}
\end{equation*}
and \prettyref{prop:pull-back} implies the following identity in $\Hom(
\propF_\emptyset(0,n)^{\Perm{n}},\propF_\Gamma(0))$:
\begin{equation*}
m\left(\int_{\pi^{-1}(\Gamma)}j\ud\mu
\right)
=n!\cdot m\left(\int_{\Gamma}j\ud\mu
\right)
\end{equation*}
and \prettyref{eq:first-gr-avg-two} is proven.
\medskip
\par
Another point of view on diagrams in $\catF_\Gamma(0)$ is the following:
any Feynman diagram is built up of vertices joined by edges; a Feynman
diagram containing $\Gamma$ as a distinguished sub-diagram is built by
joining edges and vertices to $\Gamma$. This implies
\begin{equation}\label{eq:for-the-corollary-two}
\int_{\Fey_{\Gamma}(0)}j\ud\mu=
\left(\int_{\Gamma}j\ud\mu\otimes\exp\left\{   
\int_{{\mathcal
V}}j\ud\mu
\right\}\right)\circ\int_{\mathcal
E(0,*)}j\ud\mu\,,
\end{equation}  
where $\mathcal V$ and $\mathcal E$ stand for ``vertices'' and ``edges''
respectively.
To give a formal proof of equation \prettyref{eq:for-the-corollary-two} we
need some definitions. 
The groupoid \(\catE\) is the subgroupoid
 of \(\catF\) consisting of Feynman diagrams
whose connected components are edges with two distinct endpoints.
It is immediate from the definition of $\catE$ that
\(\catE(0,2n+1)\) is empty. Moreover, any object of
\(\catE(0,2n)\) has only trivial automorphisms, so that
for any positive integer \(n\)
\begin{equation}\label{eq:pre-wick}
\int_{\mathcal
E(0,2n)}j\ud\mu=\sum_{p\in P_{2n}}
\left(\coprod_{i=1}^n {\xy\vloop-,(-0.1,.7)*{{{p^{}_{2i-1}}_{\text{out}}}}%
,(1.3,.7)*{{p^{}_{2i}}_{\text{out}}}\endxy}\right)
\end{equation}
where \(p\) ranges in the set $P_{2n}$ of all partitions
\(\{p_1,p_2\},\{p_3,p_4\},\dots\) of
\(\{1,2,\dots,2n\}\) in 2-element subsets.
For instance, for \(n=2\) we have
\begin{equation*}
   \int_{\mathcal
E(0,4)}j\ud\mu=
      {\xy\vloop-,(0.1,.7)*{{1_{\text{out}}}}%
,(1,.7)*{2_{\text{out}}},(1.7,0.5),\vloop-,(1.8,.7)*{{3_{\text{out}}}}%
,(2.7,.7)*{4_{\text{out}}}\endxy} +
        {\xy\vloop-,(0.1,.7)*{{1_{\text{out}}}}%
,(1,.7)*{3_{\text{out}}},(1.7,0.5),\vloop-,(1.8,.7)*{{2_{\text{out}}}}%
,(2.7,.7)*{4_{\text{out}}}\endxy}+
 {\xy\vloop-,(0.1,.7)*{{1_{\text{out}}}}%
,(1,.7)*{4_{\text{out}}},(1.7,0.5),\vloop-,(1.8,.7)*{{2_{\text{out}}}}%
,(2.7,.7)*{3_{\text{out}}}\endxy}
      \,.
  \end{equation*}
Note that the element defined by
equation \prettyref{eq:pre-wick} is \(\Perm{2n}\)-invariant. The groupoid
\(\catE(0,*)\) is defined as the union \(\cup_n\catE(0,n)\). If we denote 
by $\sf V$
the groupoid whose objects are ordinary vertices, then the groupoid 
whose objects
are disjoint unions of ordinary vertices is the symmetric power 
$\text{\rm Sym}(\sf
V)$, so that
\[
\int_{\text{\rm Sym}(\mathcal V)} \varphi\, \ud
\mu^{}_{\text{\rm Sym}(\mathcal V)}=
\exp \left\{\int_{\mathcal V}\varphi\,\ud\mu^{}_{\mathcal
V}\right\}
\]
by \prettyref{prop:abstract-F-vs-Z}. Then, if we define
\(\Phi_1\circ\Phi_2\) to be zero when
\(\card{\Out(\Phi_2)}\neq\card{\In(\Phi_1)}\), equation
\prettyref{eq:for-the-corollary-two} follows reasoning as in the proof 
of equation
\prettyref{eq:first-gr-avg-two}.

A variant of equation \prettyref{eq:for-the-corollary-two} is the 
following. If
$\Gamma$ is a Feynman diagram, let
\(\overline{\Gamma}\) be the subgroupoid of
\(\catF_{\Gamma}(0)\) whose objects are the Feynman diagrams that
can be obtained by joining the legs of \(\Gamma\) by
edges in all possible ways. Such diagrams will be called the
\emph{closures} of \(\Gamma\). With these notations,
\begin{equation}\label{eq:closed}
\int_{\overline{\Gamma}}j
\ud\mu
=
\int_{\Gamma}j
\ud\mu
\circ\int_{\mathcal
E(0,*)}j\ud\mu\,.
\end{equation}
Clearly, if \(\Gamma\) has an odd number of legs, then both sides
of \prettyref{eq:closed} are zero.
\medskip
\par

\everyxy={0,<2em,0em>:,(0,0),} 
\section{Feynman algebras}\label{sec:fey-algebras}
This section is concerned with the linear representations of
Feynman diagrams. More precisely, if we denote by $\Fey$ the PROP
of sets obtained from the 2-PROP $\catF$ by taking isomorphism classes of
morphisms as morphisms, then the functor \emph{free vector space} (over
the field \(\setC\)) changes $\Fey$ into a PROP $\propF$ of vector spaces:
explicitly,
$\propF(m,n)$ is the free $\setC$-vector space generated by isomorphism classes of
Feynman diagrams of type $(m,n)$.
\begin{dfn}
A Feynman algebra is an algebra over the PROP $\propF$. Equivalently,
it is a
symmetric monoidal functor
\[Z\colon
\mathscr{F}\to\catVect\,,\]
where $\catVect$ denotes the category of $\setC$-vector spaces.
The morphism \(Z\) is called \emph{graphical calculus}; the operator
\(Z(\Gamma)\) corresponding to a Feynman diagram \(\Gamma\) is called
\emph{amplitude} of the diagram\footnote{By abuse of notation we will
write $Z(\Gamma)$ for $Z([\Gamma])$.}.
\end{dfn}
Since the category of objects of \(\propF\)
is \(\catN\), which is generated by
\(1\) as a symmetric monoidal category, the image of
\(Z_{\text{Ob}}\)
will be generated by the vector space \(V=Z_{\text{Ob}}(1)\),
i.e., an \(\propF\)-algebra is actually a representation
\begin{equation*}
Z\colon \propF\to \EndOp(V)\,,
\end{equation*}
where \(\EndOp(V)\) denotes the endomorphisms PROP of \(V\).
In more colloquial terms, a Feynman algebra is the datum of a
family of morphisms
\[
Z_{m,n}:{\propF}(m,n)\to \Hom(V^{\otimes m},V^{\otimes n})
\]
respecting the braiding and compositions and tensor products of diagrams.
For instance, the fact that \(Z\) must respect the braiding forces
\begin{equation*}
Z\left(\xy*!LC\xybox{%
       \vcross~{(0,1)}{(1,1)}{(0,0)}{(1,0)},(0,-.2)*{1_{\text{in}}}%
,(1,-.2)*{2_{\text{in}}}%
,(0,1.2)*{1_{\text{out}}},(1,1.2)*{2_{\text{out}}}}\endxy\right)=
\sigma^{}_{V,V}\,
\end{equation*}
namely
\begin{equation*}
Z\left(\xy*!LC\xybox{%
       \vcross~{(0,1)}{(1,1)}{(0,0)}{(1,0)},(0,-.2)*{1_{\text{in}}}%
,(1,-.2)*{2_{\text{in}}}%
,(0,1.2)*{1_{\text{out}}},(1,1.2)*{2_{\text{out}}}}\endxy\right)
(v_1\otimes v_2)=v_2\otimes v_1\,.
\end{equation*}
Informally speaking, Feynman diagrams are
freely generated by vertices and edges, i.e., a Feynman diagram can be built by
joining together a set of vertices by means of edges in a completely arbitrary way;
so one can expect that an $\propF$-algebra is equivalent to assigning in a free way
a value to an edge and to each vertex. In fact, we have the following proposition,
which can be read as a formal statement about the forementioned freeness of 
Feynman diagrams.
\begin{prop}\label{prop:Fey-alg-stru-one}
The datum of a Feynman algebra structure on a \(\setC\)-vector space
\(V\) is the datum of
\begin{enumerate}
\item a symmetric non-degenerate bilinear pairing
\(\langle\,,\,\rangle\colon V\otimes V\to
\setC\);
\item  a family
\(\{T^\alpha_{m,n}\}_{(m,n)\in\setN\times\setN}\) of tensors
\(T^\alpha_{m,n}\colon V^{\otimes m}\to V^{\otimes n}\),
indexed by
the elements \(\alpha\) of the set \(\sf Co\) of coupon colours;
\item a family
\(\{C_n^\beta\}_{n\in\setN}\) of tensors
\(C^\beta_{n}\colon V^{\otimes n}\to \setC\),
invariant with respect
to the action of the cyclic group \(\setZ/n\setZ\) on the
inputs and indexed by the elements \(\beta\) of the set
\(\sf Cy\) of cyclic colours; we call these tensors \emph{cyclic}
tensors of the Feynman algebra;
\item a family
\(\{S_n^\gamma\}_{n\in\setN}\) of tensors
\(S^\gamma_{n}\colon V^{\otimes n}\to \setC\),
 invariant with respect
to the action of the symmetric  group \(\Perm{n}\) on the
inputs and indexed by the elements \(\gamma\) of the set
\(\sf Sy\) of cyclic colours; we call these tensors \emph{symmetric}
tensors of the Feynman algebra.
\end{enumerate}
\end{prop}
A proof of this statement can be found in
\cite{murri-fiorenza;feynman}, where techniques derived by the
Reshethikin-Turaev's graphical calculus (see
\cite{bakalov-kirillov,reshetikhin-turaev;ribbon-graphs}) are used. We
remark that in \cite{murri-fiorenza;feynman} coupon, cyclic and symmetric
vertices are treated as distinct cases; yet, the arguments used there
apply to Feynman diagrams with the three kinds of vertices occurring at
the same time.
\par

Note that no relation is required between the tensors defining the
Feynman algebra structure, nor between the tensors and the bilinear
pairing
\(\langle\,,\,\rangle\). This should be thought as an algebraic counterpart
 of the fact  that Feynman
diagrams are freely generated by the vertices
(which correspond to the
tensors) and by the edges (which correspond to the pairing). 
\par In the physicists' parlance, the tensors corresponding to
vertices are called \emph{interactions} and the dual of the pairing
corresponding to the edges is called \emph{propagator}, see for
instance \cite{QFS}.  

\par
Recall that the set of colours were split into the two subsets of ordinary
and special colours, and that an identification of ordianry colours with a subset
of special colours was given. So far, we have not implemented this datum in the
definition of Feynman algebra, not to make the definition too involved. Then, let
us complete the definition of Feynman algebra, by requiring that 
an ordinary colour and the corresponding special colour define the same tensor: 
if \({\sf v}^{\alpha}_{m,n}\) is an ordinary coupon vertex of type $(m,n)$, and 
\({\sf v}^{\bmath{\alpha}}_{m,n}\) is the corresponding special vertex, then
\[
Z(v^{\alpha}_{m,n})=T^\alpha_{m,n}=T^{\bmath{\alpha}}_{m,n}
=Z(v^{\bmath{\alpha}}_{m,n})
\]
and similarly for cyclic or symmetric vertices. Actually, we will be mostly working
with Feynman algebras whose ordinary tensors depend linearly on some complex
parameter: if \(Z\colon\propF\to\EndOp(V)\) is a Feynman algebra, we denote by
\begin{equation*}
Z_{x_*}\colon\propF\to\EndOp(V)[x_*]\,.
\end{equation*}
 the Feynman algebra obtained by
changing the ordinary tensor
\(T_{m,n}^\alpha\) (resp. \(C_{n}^\beta\) and
\(S_{n}^\gamma\) ) of \(Z\) with the tensor
\(x^{\alpha}_{m,n} T_{m,n}^\alpha\) (resp.
\(x^{\beta}_{n} C_{n}^\alpha\) and
\(x^{\gamma}_{n} S_{n}^\alpha\) ), where the
\(x_*\) are complex variables, and leaving the special tensors
unchanged. 
\par
To make \(Z_{x_*}\) a graded morphism, we put the variables
\(x^\alpha_{m,n}\) in degree
\(n+m\) and the variables \(x^\beta_{n}\) and
\(x^\gamma_{n}\) in degree \(n\).

Note that, if \({\sf v}^{\alpha}_n\) is an ordinary (coupon, cyclic or
symmetric) \(n\)-valent vertex and 
\({\sf v}^{\bmath{\alpha}}_n\) is the corresponding special vertex, then by
definition of $Z_{x_*}$ we have 
$Z_{x_*}({\sf v}^{\alpha}_n)=x^\alpha_nZ_{x_*}({\sf v}^{\bmath{\alpha}}_n)$, so that
\begin{equation*}
\frac{\del}{\del x^\alpha_n}Z_{x_*}({\sf
v}^{\alpha}_n)=Z_{x_*}({\sf
v}^{\bmath{\alpha}}_n)\,
\end{equation*}
i.e., the derivatives of the amplitudes of ordinary vertices with
respect to the parameters \(x_*\) can be written as amplitudes of the
corresponding special vertices. In particular, this implies that the
derivatives of the amplitude of a Feynman diagram \(\Gamma\) with
respect to the parameters \(x_*\) can be written as sums of copies of
diagrams obtained by changing some ordinary vertex of \(\Gamma\) with
the corresponding special vertex, which justify the convention of having the
ordinary colours identified with a subset of the special ones.
\medskip

If $\Gamma$ is an object of $\catF(n)$ then, due to
the lack of an ordering on the legs,
\(\Gamma\) does not define a linear operator via the graphical
calculus
\(Z_{x_*}\). Anyway, two pre-images of
\(\Gamma\) in \(\catF(n,0)\) can only differ by the action of an
element of the symmetric group \(\Perm{n}\), so that a linear operator
$Z_{x_*}(\Gamma)$ on the subspace of \(\Perm{n}\)-invariant vectors of
 \(V\tp{n}\) is well
defined: 
\begin{equation*}
Z_{x_*}(\Gamma)=Z_{x_*}(\hat\Gamma)
\bigr\vert_{(V\tp{n})^{\Perm{n}}}
\colon(V\tp{n})^{\Perm{n}}\to\setC\,,
\end{equation*}
where \(\hat\Gamma\) is any Feynman diagram in the fibre of
\(\catF(n,0)\to\catF(n)\) over \(\Gamma\).
Hence we obtain the \emph{polynomial function} associated to the
 diagram
\(\Gamma\):
\begin{equation*}
P_{\Gamma}\colon v\mapsto Z_{x_*}^{}(\Gamma)(v^{\otimes n})
\end{equation*}
The association $\Gamma\mapsto P_\Gamma^{}$ extends to a linear
map
\begin{equation*}
P\colon\propF(n)\to\{\text {degree $n$ homogeneous polynomials on }V\}
\end{equation*}
where $\propF(n)$ denotes the $\setC$-vector space generated by
isomorphism classes of objects in $\catF(n)$. We remark that sometimes the
term ``amplitude'' is used to denote the
polynomial function $P_\Gamma$ rather than the linear operator 
$Z_{x_*}(\Gamma)$.

\everyxy={0,<2em,0em>:,(0,0.5),} 

\section{Expectation values}

Let now $Z_{x_*}\colon\catF\to\End(V)[x_*]$ be a Feynman algebra as
defined in \prettyref{sec:fey-algebras}. Since 
the morphism
$Z_{x_*}$ is graded by construction, it is integrable as a morphism
\begin{equation*}
Z_{x_*}\colon\catF\to\End(V)[[x_*]]
\end{equation*}
by \prettyref{lemma:filtrated}.
 The \emph{expectation value} of
\(\Gamma\) is the element of \(\setC[x_*]\) defined by
\begin{equation*}
\langle\!\langle\Gamma\rangle\!\rangle
=\int_{\overline{\Gamma}}
Z_{x_*}\ud\mu\,,
\end{equation*}
or, in the more familiar ``sums'' notation,
\begin{equation*}
\langle\!\langle\Gamma\rangle\!\rangle
=\sum_{[\Phi]\in\overline{\Gamma}}\frac{Z_{x_*}(\Phi)}
{\card{\Aut\Phi}}.
\end{equation*}
The \emph{expectation value with
potential}\footnote{This terminology will be explained in
\prettyref{sec:expansion-integrals}.} of
\(\Gamma\) is the formal series in the variables \(x_*\) defined
by
\begin{equation*}
\langle\!\langle\Gamma\rangle\!\rangle^{}_{x_*}:=
\int_{\Fey_{\Gamma}(0)}
Z_{x_*}\ud\mu\,,
\end{equation*}
or, in the more familiar ``sums'' notation,
\begin{equation*}
\langle\!\langle\Gamma\rangle\!\rangle^{}_{x_*}:=
\sum_{[\Phi]\in\Fey_{\Gamma}(0)}
\frac{Z_{x_*}(\Phi)}{\card{\Aut\Phi}}\,.
\end{equation*}
\par
A useful relation among expectation values is the following. Assume
that \(\Gamma\) and \(\Phi_1,\dots\Phi_k\) are
objects of \(\catF(n)\), for some fixed \(n\), and that,
for suitable polynomials \(a_k(x_*)\)
\begin{equation}\label{eq:give-same-expect}
\int_{\Gamma}Z_{x_*}^{}\ud\mu=\sum_k
a_k(x_*)\int_{\Phi_k}Z_{x_*}^{}\ud\mu
\end{equation}
as linear operators on \(\left(V\tp{n}\right)^{\Perm{n}}\). Then, by
equations
\prettyref{eq:for-the-corollary-two} and \prettyref{eq:closed}, it 
immediately
follows
\begin{equation*}
\langle\!\langle \Gamma \rangle\!\rangle=
\sum_k
a_k(x_*)\langle\!\langle
\Phi_k\rangle\!\rangle
\quad\text{ and }\quad
\langle\!\langle \Gamma\rangle\!\rangle_{x_*}=\sum_k
a_k(x_*)\langle\!\langle
\Phi_k\rangle\!\rangle_{x_*}\,.
\end{equation*}
In the more familiar ``sums'' notation, equation
\prettyref{eq:give-same-expect} above is the following identity among
polynomials on \(V\):
\begin{equation*}
\frac{P_{\Gamma}(v)}{\card{\Aut\Gamma}}=\sum_k
a_k(x_*)\frac{P_{{\Phi}_k}(v)}{\card{\Aut{\Phi}_k}}\,.
\end{equation*}

\section{Partition functions and free energy}\label{sec:pffe}

We have seen that the graphical calculus \(Z_{x_*}\) is integrable on
\(\catF_{\Gamma}\), for any Feynman diagram \(\Gamma\). The
\emph{partition function}  of the Feynman algebra
\(Z_{x_*}\colon \catF\to \EndOp(V)[x_*]\) is the formal series
defined as the integral of \(Z_{x_*}\) on the groupoid
\(\catF_{\emptyset}(0)\) of all Feynman diagrams with no legs
and only ordinary vertices, 
namely,
\begin{equation}
\boxed{
Z(x_*):=
\int_{\Fey_{\emptyset}(0)}Z_{x_*}\ud\mu=\langle\!\langle\emptyset
\rangle\!\rangle^{}_{x_*}}
\end{equation}
The \emph{free energy} of the Feynman algebra \(Z_{x_*}\) is defined as the
integral of \(Z_{x_*}\) on the subgroupoid \(\catF_{\emptyset,\text{\rm 
conn.}}(0)\)
of \(\catF_{\emptyset}(0)\) consisting of connected Feynman diagrams 
with only
ordinary vertices: 
\begin{equation}
 F(x_*):=
\int_{\Fey_{\emptyset,\text{\rm
conn.}}(0)}Z_{x_*}\ud\mu\,,
\end{equation}
In the more familiar ``sums'' notation
\begin{align*}
Z(x_*)&=
\sum_{[\Phi]\in\Fey_{\emptyset}(0)}\frac{Z_{x_*}(\Phi)}
{\card{\Aut\Phi}}
\\ F(x_*)&=
\sum_{[\Phi]\in\Fey_{\emptyset,\text{\rm
conn.}}(0)}\frac{Z_{x_*}(\Phi)}{\card{\Aut\Phi}}
\end{align*}

The following well known relation connects the partition function with
the free energy:
\begin{equation}\label{eq:FvsZ}
\boxed{Z(x_*)=\exp\bigl\{F(x_*)\bigr\}}
\end{equation}
To prove  it, just
observe that there is an isomorphism of measure
spaces
\begin{equation*}
\Fey_\emptyset(0)\simeq
\text{\rm{Sym}}\bigl(\Fey_{\emptyset,\text{\rm conn.}}(0)\bigr)
\end{equation*}
and apply
\prettyref{prop:abstract-F-vs-Z}.

We now show how the derivatives of the partition function 
$Z(x_*)=\langle\!\langle
\emptyset \rangle\!\rangle^{}_{x_*}$ are related to the expectation 
values of
disjoint unions of special vertices. Let
\(\alpha\) be an ordinary colour for a vertex of
\(\catF\) and let
\({\sf v}_{}^\alpha\) and
\({\sf v}_{}^{\bmath{\alpha}}\) be the corresponding
ordinary
and special vertices. By definition
of $Z_{x_*}^{}$, we have
\begin{equation*}
Z_{x_*}^{}({\sf v}^\alpha_{})=x^{}_\alpha\cdot
Z_{x_*}({\sf v}^{\bmath{\alpha}}_{})\,,
\end{equation*}
with \(Z_{x_*}({\sf v}^{\bmath{\alpha}}_{})\) which is actually
independent of the variables \(x_*\). So, if we apply the differential
operator \(\del/\del x_{}^\alpha\) to the partition function \(Z(x_*)\), we
find
\begin{equation*}
\frac{\del Z(x_*)}{\del
x_{}^\alpha}=
\frac{\del}{\del x_{}^\alpha}\left(\left(
\exp\left\{
\int_{{\mathcal
V}}Z_{x_*}\ud\mu\right\}\right)\circ\int_{\mathcal
E(0,*)}Z_{x_*}\ud\mu\right)
\end{equation*}
If \(\sf e\) is an edge, the graphical calculus \(Z_{x*}^{}(\sf e)\)
is actually independent of \(x_*^{}\), and the right-hand side of the
above equation equals
\begin{align*}
&\left(
\left(
\int_{{\sf v}_{}^{\bmath{\alpha}}}Z_{x_*}\ud\mu\right)\otimes
\left(\exp\left\{
\int_{{\mathcal
V}}Z_{x_*}^{}\ud\mu
\right\}
\right)\right)
\circ\int_{\mathcal
E(0,*)}Z_{x_*}\ud\mu
\\
&\qquad=
Z_{x_*}^{}\left(\left(
\int_{{\sf v}_{}^{\bmath{\alpha}}}j\ud\mu\otimes
\exp\left\{
\int_{{\mathcal
V}}j\ud\mu
\right\}
\right)
\circ\int_{\mathcal
E(0,*)}j\ud\mu\right)\\
&\qquad=Z_{x_*}^{}\int_{{\mathcal F}_{{\sf
v}_{}^{\bmath{\alpha}}}(0)}j\ud\mu\,,
\end{align*}
by equation \prettyref{eq:for-the-corollary-two}.
But this is just the expectation value with potential of the vertex
\({\sf v}_{}^{\bmath{\alpha}}\). So, summing up, we have
proved the formula
\begin{equation}\label{eq:derivative-vacuum}
\boxed{
\frac{\del \langle\!\langle\emptyset\rangle\!\rangle^{}_{x_*}}{\del
x_{}^\alpha}=\langle\!\langle{\sf
v}_{}^{\bmath{\alpha}}\rangle\!\rangle^{}_{x_*}
}
\end{equation}
which holds for any ordinary colour \(\alpha\) for the vertices of
\(\catF\).
More generally,
\begin{equation}\label{eq:derivative-vacuum-two}
\boxed{
\frac{\del^e_{} \langle\!\langle\emptyset\rangle\!\rangle^{}_{x_*}}{({\del
x_{}^{\alpha_1}})^{e_1}\cdots({\del
x_{}^{\alpha_k}})^{e_k}}=(e_1!\cdots e_k!)\langle\!\langle{\sf
v}_{}^{\bmath{\alpha}_1} \otimes\cdots\otimes {\sf
v}_{}^{\bmath{\alpha}_k}\rangle\!\rangle^{}_{x_*}
}
\end{equation}
for any choice of ordinary colours  \(\alpha_1^{},\dots,\alpha_k^{}\) for
the vertices of
\(\catF\).

We end this section with a digression on $\Gamma$-reduced diagrams.
An 
object \(\Phi\)
of \(\catF_{\Gamma}\) is called
\(\Gamma\)-\emph{reduced} if no connected component of \(\Phi\) has
empty intersection with \(\Gamma\). In particular, the only
\(\emptyset\)-reduced Feynman diagram is the empty Feynman diagram.
We denote the subgroupoid of \(\Gamma\)-reduced Feynman diagrams
by the symbol
\(\catF_{\Gamma\text{\rm -red.}}\).
The
isomorphism of measure spaces
\begin{equation*}
\Fey_{\Gamma}(0)\simeq
\Fey_{\Gamma\text{\rm -red.}}(0)\times\Fey_\emptyset(0)
\end{equation*}
gives
\begin{equation}\label{eq:gamma-reduced}
\int_{\Fey_{\Gamma\text{\rm -red.}}(0)}Z_{x_*}\ud\mu
=
\frac{1}{Z(x_*)}\int_{\Fey_{\Gamma}(0)}
Z_{x_*}\ud\mu  
\,,
\end{equation}
i.e, in the ``sums'' notation,
\begin{equation*}
\sum_{[\Phi]\in{\Fey_{\Gamma\text{\rm
-red.}}(0)}}\frac{Z_{x_*}(\Phi)}{\card{\Aut\Phi}}
=
\frac{1}{Z(x_*)}\sum_{[\Phi]\in{\Fey_{\Gamma}(0)}}
\frac{Z_{x_*}(\Phi)}{\card{\Aut\Phi}}\,.
\end{equation*}

\section{Gaussian integrals}              

Let \(V\) be a finite dimensional real Hilbert space,
with inner product \(\inner{-}{-}\). If \(\{e_i\}\) is a basis of
\(V\), we denote the \emph{coordinate} maps relative to this basis as
\(e^i:V\to\setR\), and write \(v^i=e^i(v)\), for any vector
\(v\) of \(V\). Via the inner product of \(V\), we can
identify the functionals
\(\{e^i\}\) with vectors of \(V\) that we will denote by
the same symbols. The vectors
\(\{e^i\}\) are a basis for \(V\), called the dual basis with
respect to \(\{e_i\}\).
 The matrix associated to
\(\inner{-}{-}\) with respect to the basis
\(\{e_i\}\) is the matrix \((g_{ij}^{})\) defined  by
\begin{equation*}
  g^{}_{ij} := \inner{e_i}{e_j}.
\end{equation*}
As customary, we set \(g^{ij} := (g^{-1})_{ij} = \inner{e^i}{e^j}\).

Let now \(\ud v\) be a (non trivial) translation invariant measure on
\(V\). The function \(\E^{-\onehalf\inner{v}{v}}\) is positive and
integrable with respect to \(\ud v\).
  The probability measure on \(V\) defined by
  \begin{equation*}
\ud\mu(v)=\frac{\displaystyle{\E^{-\frac{1}{2}
\inner{v}{v}}\ud
v}}{\displaystyle{\int_{V}
\E^{-\frac{1}{2}
\inner{v}{v}}\ud
v}}\,.
\end{equation*}
  is called the \emph{Gaussian measure} on \(V\).
Since a non-trivial translation invariant measure on \(V\) is unique up
to a scalar factor, \(\ud\mu\) is actually independent of the chosen
\(\ud v\).

The inner product \(\inner{-}{-}\) extends uniquely to a symmetric
\(\setC\)-bilinear pairing on the complex vector space
\(\cplx{V} := V \otimes
\setC\); this pairing, which we shall denote by
the same symbol \(\inner{-}{-}\), is clearly
non-degenerate. Identify
\(V\) with the subspace
\(V
\otimes
\{1\}\)  of real vectors in
\(\cplx{V}\). Polynomial functions on \(\cplx{V}\) are
integrable with respect to the Gaussian measure; for any
polynomial function
\(f\colon V_\setC\to\setC\) we set
\begin{equation*}
\gint{f}=\int_Vf(v)\ud\mu(v)\,.
\end{equation*}
The complex number \(\gint{f}\) is called the \emph{average}
of \(f\) with respect to the Gaussian measure.

Since the vectors \(\{ e_i \}\) are a basis for the complex
vector space \(\cplx{V}\), the vectors
\(\{e_{i_1}\otimes\cdots\otimes e_{i_n}\}\) are a basis for
the vector space \(\cplx{V}\tp{n}\): any elment \(v_{[n]}\)
of
\(\cplx{V}\tp{n}\) can be uniquely written as
\begin{equation*}
v_{[n]}=\sum_{i_1,\dots,i_n}v_{[n]}^{i_1,\dots,i_n}e_{i_1}\otimes\cdots\otimes
e_{i_n}\,.
\end{equation*}
For any \(v\in \cplx{V}\), the vector \(v\tp{n}\) is an
element of \(\cplx{V}\tp{n}\), and the following identity
holds:
\begin{equation*}
v\tp{n}=\sum_{i_1,\dots,i_n}v^{i_1}\cdots
v^{i_n}e_{i_1}\otimes\cdots\otimes e_{i_n}\,.
\end{equation*}
The functions \(v\mapsto v^{i_1}\cdots
v^{i_n}\) are polynomials on \(\cplx{V}\) and
we can define the average of \(v\tp{n}\) as
\begin{equation*}
\gint{v\tp{n}}:=\sum_{i_1,\dots,i_n}
\gint{v^{i_1}\cdots v^{i_n}}e_{i_1}\otimes\cdots\otimes
e_{i_n}\,;
\end{equation*}
it is clearly independent of the basis \(\{e_i\}\) chosen.

The golden bridge between Gaussian integrals and Feynman diagrams is the
following Lemma, due to Gian Carlo Wick. In its original formulation it
is stated in terms of momenta of the Gaussian measure, i.e., averages of
monomials in the coordinates \(v^i\); see for instance
\cite{bessis-itzykson-zuber;graphical-enumeration}. Here, using the
notion of average of  \(v\tp{n}\), we recast it in
a coordinate-free way, which is more suitable for a reinterpretation
through the graphical formalism of the previous sections. Recently,
Robert
Oeckl has proven that a Wick-type lemma holds in the wider context of
general braided tensor categories, see
\cite{oeckl;braided-qft}.
\begin{lemma}[Wick]
  Tensor powers of vectors are integrable with
  respect to the
Gaussian measure \(\ud\mu\) and:
  \begin{align}
  \label{eq:W1}
  {\avg{v\tp{ 2n+1}}} &= 0,
  \\
  \label{eq:W2}
  {\avg{v\tp{ 2}}} &= \sum_{i,j} g^{ij} e_i
  \otimes e_j,
  \\
  \label{eq:W3}
  {\avg{v\tp{2n}}} &= \sum\limits_{i_1, \dots,
    i_{2n}} \sum\limits_{s\in P} g^{i_{s_1} i_{s_2}} \cdots
  g^{i_{s_{2n-1}} i_{s_{2n}}} e_{i_1} \otimes \cdots
  \otimes e_{i_{2n}},
\end{align}
  where \(P\) denotes the set of all distinct pairings of the set of
  indices \(\{i_1,\dots,i_{2n}\}\), i.e., over the set of all
  partitions \(\{\{i_{s_1},i_{s_2}\},\{i_{s_3},i_{s_4}\},\dots\}\) of
  \(\{i_1,i_2,\dots,i_{2n}\}\) into 2-element subsets.
\end{lemma}

\section{Feynman diagrams expansion of Gaussian
integrals}\label{sec:expansion-integrals}           

\everyxy={0,<2em,0em>:,(0,.5),} 

We now show how a Gaussian integral can be
expanded into a sum of Feynman diagrams, to be evaluated
according to the rules of graphical calculus. Historically,
Gaussian integrals are the context where Richard Feynman
originally introduced the diagrams that nowadays bear his name. The key 
point will
be a graphical interpretation of Wick's lemma.

Let \(Z_{x_*}\colon \catF\to \EndOp(\cplx{V})[x_*]\) be a Feynman algebra
compatible with the bilinear pairing \(\inner{-}{-}\), i.e., such
that
\begin{equation*}
{\xy
\vloop~{(0,0.5)}{(1,0.5)}{(0,-.5)}{(1,-.5)},(0,-.7)*{
x} ,(1,-.7)*{y}
\endxy}=\inner{x}{y},\qquad\forall x,y\in \cplx{V}\,.
\end{equation*}
Since the right-hand side of (\ref{eq:W2}) in Wick's lemma is the 
co-pairing
relative to the pairing  \(\inner{-}{-}\) on \(\cplx{V}\), we can rewrite
\prettyref{eq:W2} as
\begin{equation*}
  {\avg{v\tp{2}}}=
Z_{x_*}\left({\xy\vloop-,(0.1,.7)*{{1_{\text{out}}}}%
,(1,.7)*{2_{\text{out}}}\endxy}\right)=\int_{{\mathcal
E}(0,2)}Z_{x_*}\ud\mu\,
.
\end{equation*}
Also \prettyref{eq:W3} can be expressed as an integral over Feynman 
diagrams: by
\prettyref{eq:pre-wick}, we have the following recasting of
 Wick's lemma in terms of integrals on the groupoid of edges.
\begin{lemma}[Wick's lemma via graphical calculus]
  Tensor powers of vectors are integrable with
  respect to the
Gaussian measure \(\ud\mu\) and:
\begin{equation}
  \label{eq:avg-to-casimir}
  {\avg{v\tp{m}}}
=\int_{{\mathcal
E}(0,m)}Z_{x_*}\ud\mu\,
.
\end{equation}
\end{lemma}

We now introduce the potential of a Feynman algebra. Recall that
\(\catV\) denotes the groupoid whose objects are ordinary vertices with no
numbering on the legs.  The \emph{potential} of the  Feynman algebra
\(Z_{x_*}^{}\colon\catF\to\EndOp(V)[x_*]\) is defined as the formal series
\begin{equation}\label{eq:potential}
S(x_*):=
\int_{\mathcal
V} P\ud\mu\,,
\end{equation}
where $P$ denotes the polynomial function associated to a diagram.
In more familiar ``sums'' notations, equation
\prettyref{eq:potential} above reads
\begin{align*}
S(x_*;v)&:=
\sum_{[{\sf v}]\in {\mathcal
V}}\frac{P_{\sf
v}(v)}{\card{\Aut{\sf v}}}\\
&=\sum_{{\substack{m,n\in \setN\\ \alpha\in{\sf Co}}}}x_{m,n}^\alpha
\bigl\langle
T_{m,n}^\alpha(v^{\otimes m}),v^{\otimes n}
\bigr\rangle+\sum_{{\substack{n\in \setN\\ \beta\in{\sf
Cy}}}}x_{n}^\beta
\frac{C_{n}^\beta(v^{\otimes n})}{n} +\sum_{{\substack{n\in\setN\\
\gamma\in{\sf
Sy}}}}x_{n}^\gamma
\frac{S_{n}^\gamma(v^{\otimes n})}{n!}
\end{align*}

If \(f\) is a polynomial function on \(V\), the average of
\(f\) with potential \(S(x_*)\) is the formal series in the
variables \(x_*\) defined by
\begin{equation}
{\gint{f}}_{x_*}=\gint{f\cdot \E^{S(x_*)}}
\end{equation}

If \(\Gamma\) is a Feynman diagram, the function
\(P_{\Gamma}\) is a  polynomial on
\(V\). Therefore we can consider the average of \(P_{\Gamma}\)
with respect to the Gaussian measure. The following result shows that
our
notations are consistent.
\begin{thmsec}[Feynman-Reshetikhin-Turaev]\label{thm:FRT}
 For any Feynman diagram \(\Gamma\)
 the following equations hold:
 \begin{equation}
    \label{eq:FRT0}
  {\langle\!\langle\Gamma\rangle\!\rangle}=\left\langle
\frac{P_{\Gamma}}{\card{\Aut\Gamma}}\right\rangle=
  \int_V\frac{P_{\Gamma}(v)}{\card{\Aut\Gamma}}
\ud\mu^{}_V(v)\,;
  \end{equation}  
\begin{equation}
    \label{eq:FRT1}
  {\langle\!\langle\Gamma\rangle\!\rangle}_{x_*}=\left\langle
\frac{P_{\Gamma}}{\card{\Aut\Gamma}}\right\rangle^{}_{x_*}=
  \int_V\frac{P_{\Gamma}(v)}{\card{\Aut\Gamma}}
\E^{S(x_*;v)}
\ud\mu^{}_V(v)\,.
  \end{equation}
\end{thmsec}
\begin{proof}
We prove only \prettyref{eq:FRT0}, equation
\prettyref{eq:FRT1} being completely analogous. Let \(n\) be the number
of legs of \(\Gamma\); by linearity
\begin{equation*}
\frac{1}{\card{\Aut{\Gamma}}}\langle
P_{\Gamma}\rangle=
\frac{Z_{x_*}(\Gamma)}{\card{\Aut{\Gamma}}}{\gint{v\tp{n}}}
\end{equation*}
If \(n\) is odd, the right-hand side of the above equation is zero by
\prettyref{eq:avg-to-casimir}. On the other hand, a Feynman diagram with
an odd
number of legs cannot be closed by joining its endpoints by edges,
since a disjoint union of edges always has an even number of
endpoints; so equation
\prettyref{eq:FRT0} is verified for odd \(n\). For \(n=2m\) we
find, again by
\prettyref{eq:avg-to-casimir} and by \prettyref{eq:closed},
\begin{align*}
\frac{Z_{x_*}(\Gamma)}{\card{\Aut{\Gamma}}}
{\gint{v\tp{n}}}&=\int_{\Gamma}Z_{x_*}\ud\mu\circ
\left(
\int_{{\mathcal E}(0,2m)}
Z_{x_*}\ud\mu\right)\\
&=Z_{x_*}\left(
\int_{\Gamma}j\ud\mu\circ
\int_{{\mathcal E}(0,2m)}
j\ud\mu\right)
\\
&=Z_{x_*}
\int_{\overline{\Gamma}}j\ud\mu
={\langle\!\gint{\Gamma}}\!\rangle
\end{align*}
\end{proof}

As a particular case of \prettyref{eq:FRT1}, corresponding to 
$\Gamma=\emptyset$, we
have
\[
Z(x_*)=\int_V
\E^{S(x_*;v)}
\ud\mu^{}_V(v)\,.
\]
\medskip

 We will conclude with three examples of equation \prettyref{eq:FRT0}
involving a coupon, a cyclic or a symmetric vertex. 
Equation
\prettyref{eq:avg-to-casimir}, for
\(n=2\) gives

\hskip -1 em\includegraphics{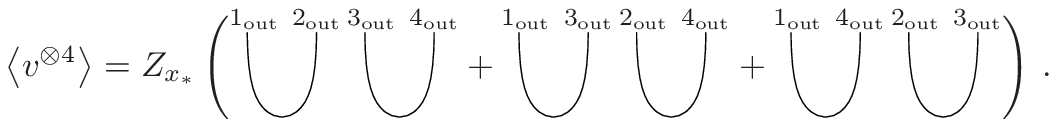}

\noindent
Consider now a \(4\)-valent coupon
vertex decorated by the colour \(\alpha\). By linearity,

\medskip

\centerline{\includegraphics{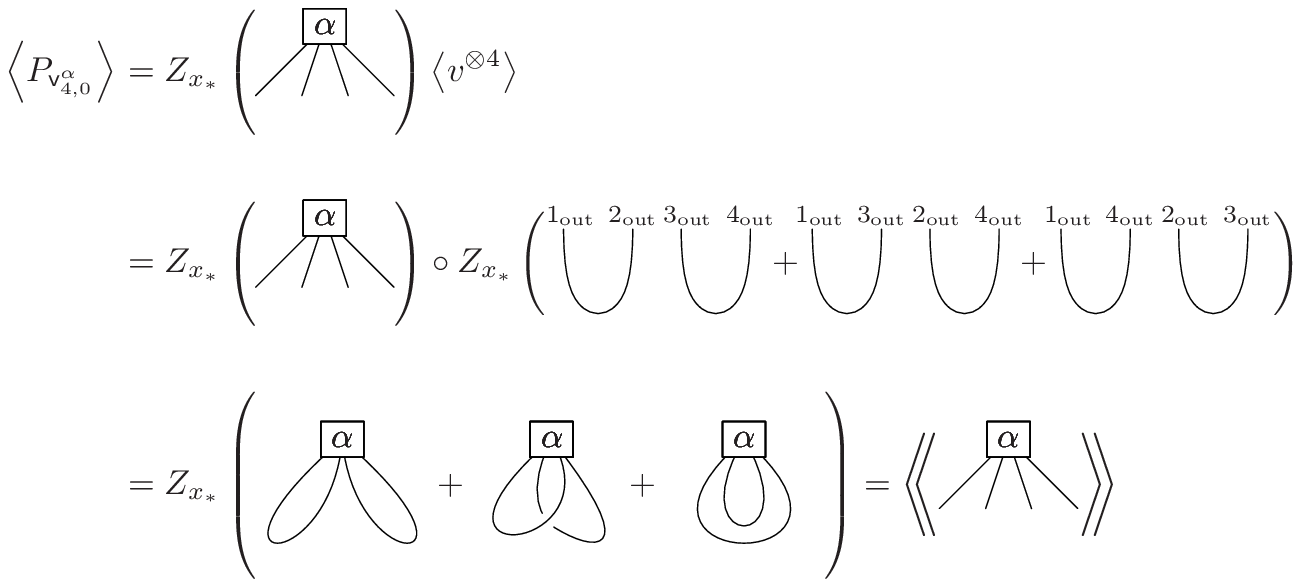}}

\noindent
If we  consider a \(4\)-valent cyclic
vertex decorated by the colour $\beta$ instead, we  have

\centerline{\hskip 1 em\includegraphics{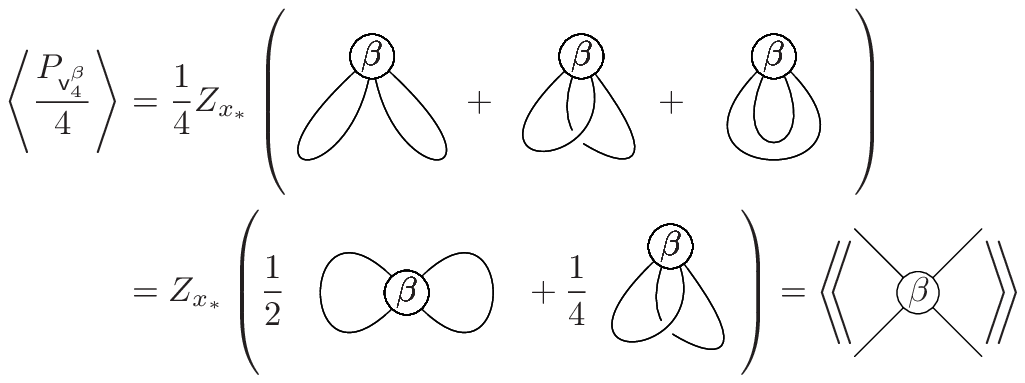}}

\noindent
Finally, if we consider a  symmetric
\(4\)-valent vertex decorated by the colour $\gamma$, we find

\centerline{\hskip 1 em\includegraphics{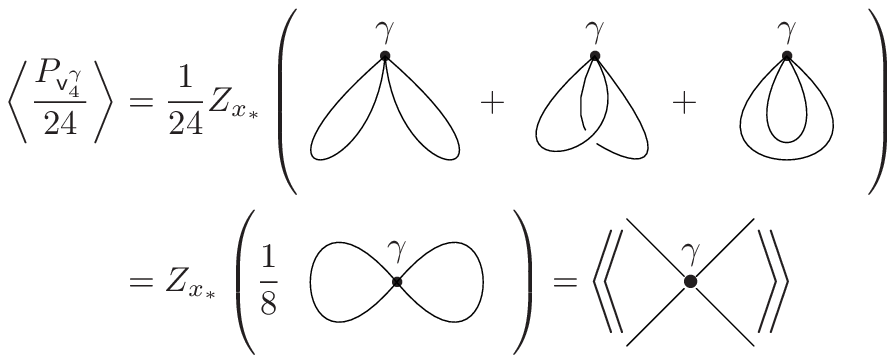}}

\section*{Appendix: Working with coordinates}

\everyxy={0,<2em,0em>:,(0,0.3),} 

In the main body of the paper we have only worked with
coordinate-free formulas; since the pairing \(v\otimes
w\mapsto(v,w)\) of a Feynman algebra is a nondegenerate symmetric pairing on
a $\setC$-vector space $V$, it 
 admits orthonormal bases
\(\{e_i\}_{i=1\dots N}\), and it is possible to write the coordinate version of
these formulas with respect to an orthonormal basis \(\{e_i\}\) in terms of Feynman
diagrams with coloured edges. In particular, when the Feynman algebra
is related to the Feynman diagram expansion of a Gaussian integral, the vector space
$V$ is the the complexification of a real Hilbert
space $V^{}_{\setR}$ and one can choose $\{e_i\}$ to be the complexification of an
othonormal basis of $V_{\setR}^{}$.

By definition, a Feynman diagram with edges coloured by the symbols
\(\{1,\dots,N\}\) is a pair \((\Gamma,\eta)\), where \(\Gamma\) is a
Feynman diagram and \(\eta\) is a map \(\eta\colon{\rm Edges}(\Gamma) \to
\{1,\dots,N\}\). We represent an edge coloured by the symbol ``\(i\)'' by
writing an ``\(i\)'' near it. The groupoid of Feynman diagrams
with edges coloured by \(\{1,\dots,N\}\) is denoted by the symbol
\({}_{\{1,\dots,N\}}\catF\); clearly, automorphisms of Feynman diagrams
with coloured edges are required to preserve the colouring. It is
immediate to check that the  ``forget the colouring on the edges'' map is a
covering
\begin{equation*}
\pi\colon {}_{\{1,\dots,N\}}\catF\to \catF
\end{equation*}

 We now have to define the amplitudes
\(Z_{x_*}^{}(\Gamma,\eta)\). Let \(\pi_i\colon V\to V\) be the orthogonal
projection on the subspace spanned by \(e_i\). We make the graphical
assignment
\begin{equation*}
{\xy*!LC\xybox{
        (0,-.7);(0,0.7)**\dir{-}
        ,(0,-.9)*{1_{\text{\rm in}}}%
        ,(0,0.9)*{1_{\text{\rm out}}}%
        ,(.15,-.3)*{i}
        ,(1.5,-.1)*{\displaystyle{\mapsto \pi_i\,.}}
        }\endxy}
\end{equation*}
Since the basis \(\{e_i\}\) is orthonormal with respect to the pairing
\((-,-)\), this graphical assignment is consistent with the other
rules of graphical calculus and so a well-defined amplitude
$Z_{x_*}$ is induced on Feynman diagrams with coloured edges.
 
The equation \(\Id_V=\oplus_{i=1}^N\pi_i\) is translated in graphical
terms into
\begin{equation*}
{\xy*!LC\xybox{
        (0,-.7);(0,0.7)**\dir{-}
        ,(0,-.9)*{1_{\text{\rm in}}}%
        ,(0,0.9)*{1_{\text{\rm out}}}%
 }\endxy}       
\begin{matrix}=&\displaystyle{\bigoplus_{i=1}^N}\\{}\\\end{matrix}
{\xy*!LC\xybox{
        (0,-.7);(0,0.7)**\dir{-}
        ,(0,-.9)*{1_{\text{\rm in}}}%
        ,(0,0.9)*{1_{\text{\rm out}}}%
        ,(.15,-.3)*{i}
}\endxy}
\end{equation*}
so that the amplitude of a Feynman diagram is expanded into a sum of
amplitudes of Feynman diagrams with coloured edges:
\begin{equation}\label{eq:sum-sum-col}
Z_{x_*}(\Gamma)=\sum_{\eta\colon{\rm Edges}(\Gamma)\to\{1,\dots N\}}
Z_{x_*}(\Gamma,\eta)
\end{equation}
But this is just a push-forward formula:
\begin{equation*}
Z_{x_*}=\pi_*Z_{x_*}^{}
\end{equation*}
where \(\pi\) is the ``forget the colouring on the edges'' map. So we can
apply Fubini's theorem (\prettyref{prop:fubini}) and find
\begin{equation}\label{eq:fubini-fey}
\int_{\mathcal F}Z_{x_*}^{}\ud\mu_{\mathcal F}= \int_{{}_{\{1,\dots
N\}}{\mathcal F}}Z_{x_*}^{}\ud\mu_{{}_{\{1,\dots
N\}}{\mathcal F}}
\end{equation}
that is, summing over Feynman diagrams is the same thing as summing over
Feynman diagrams with coloured edges (obviously, with the right weights).

 A classical example of
\prettyref{eq:fubini-fey} is the following. Let \(\phi\) be an analytic
function defined on a neighborhood of $0$ in \(V\). The derivatives
\(D^n\phi\) define
a family of symmetric tensors on \(V\). We make the graphical assignment
\everyxy={0,<2em,0em>:,(0,0.5),} 
\[
{\xy*!LC\xybox{
        \blackvertex{7}\loose1\loose2\loose3\loose4\missing5\loose6%
\loose7%
,(1.55,.2)*{1_{\text{\rm in}}}
,(1.3,1.3)*{2_{\text{\rm in}}}
,(-1.55,.2)*{5_{\text{\rm in}}}
,(-1.3,1.3)*{4_{\text{\rm in}}}
,(0,2)*{3_{\text{\rm in}}}
,(0.9,-.8)*{n_{\text{\rm in}}}
        }\endxy}\quad\begin{matrix}\mapsto
\displaystyle{D^n\phi\bigr\rvert_0}
\\{}\\
\end{matrix}
\]
\everyxy={0,<2em,0em>:,(0,0.3),} 
and call this an $n$-valent black vertex.
Any vector $v\in V$ can be seen as a morphism $v\colon \setC\to V$; we
make the graphical assignment
\[
{\xy*!LC\xybox{
        ,(0,0)*+[F]{v};(0,1.2)**\dir{-}
,(0,1.4)*{1_{\text{\rm out}}}
}        \endxy}\quad\begin{matrix}\mapsto \displaystyle{v}
\\{}\\
\end{matrix}
\]
and call this a  $v$-vertex. Both types of vertices will be
considered as special. Finally, for a fixed $v\in V$, we denote by
${\sf Stars}(v)$ the groupoid of Feynman diagrams whose objects are
the diagrams with no legs and exactly one
black vertex such that all the edges stemming from the black vertex end
in some $v$-vertex. 
If $v$ is a vector in the domain of $\phi$, the Taylor formula for
$\phi$ can then be written as 
a sum over
Feynman diagrams:
\begin{equation}\label{eq:taylor-one}
\phi(v)=\int_{{\text{\sl Stars}}(v)}Z_{x_*}\ud\mu=
\int_{{}_{\{1,\dots N\}}{\text{\sl Stars}}(v)}Z_{x_*}\ud\mu\,.
\end{equation}
Written out explicitly, the above equation is just the well-known identity
\begin{multline}\label{eq:taylor-two}
\phi(v)=\sum_{n=0}^\infty \frac{1}{n!}D^n\phi\bigr\rvert_0(v^{\otimes
n})=\\
=
\sum_{d_1,\dots d_N=0}^\infty\frac{1}{d_1!\cdots d_N!}\cdot
\frac{\partial^{d_1+\cdots d_N}\phi}{(\partial {v^1})^{d_1}\cdots
(\partial {v^N})^{d_N}}\biggr\rvert_0 \cdot {(v^1)}^{d_1}\cdots
{(v^N)}^{d_N}\,,
\end{multline}
where \((v^1,\dots v^N)\) are the coordinates of the vector $v$ with
respect to the basis \(\{e_i\}\).


\newcommand{\etalchar}[1]{$^{#1}$}

\end{document}